%% file: main.tex
\documentclass[twocolumn]{autart}    % Enable this line and 

\usepackage{graphics} 
\usepackage{times} 
\usepackage{cite}
\usepackage{array}

\usepackage{amssymb} 
\usepackage{ulem}

\newcommand{\mR}{\mathbb{R}}
\newcommand{\mE}{\mathbb{E}}
\newcommand{\cN}{\mathcal{N}}

\newcommand{\mP}{\mathbb{P}}
\newcommand{\mQ}{\mathbb{Q}}
\newcommand{\tr}{\operatorname{Tr}}
\renewcommand{\arraystretch}{1.0}

\newtheorem{theorem}{Theorem}
\newtheorem{lemma}{Lemma}
\newtheorem{proposition}{Proposition}

\newtheorem{problem}{Problem}
\newtheorem{remark}{Remark}
\newenvironment{proof}[1][Proof]{\noindent\textbf{#1:} }{\hfill$\square$}

\usepackage{subcaption} 
\usepackage{color}

\usepackage{mathtools}

\begin{document}

\begin{frontmatter}

\title{Optimal Covariance Steering of Linear Stochastic Systems with Hybrid Transitions} 

\thanks[footnoteinfo]{Financial support from NSF under grants 1942523, 2206576 are greatly acknowledged.}

\author[Gatech]{Hongzhe Yu}\ead{hyu419@gatech.edu},    % Add the 
\author[CMU]{Diana Frias Franco}\ead{dfriasfr@andrew.cmu.edu},               % e-mail address 
\author[CMU]{Aaron M. Johnson}\ead{amj1@andrew.cmu.edu},  % (ead) as shown
\author[Gatech]{Yongxin Chen}\ead{yongchen@gatech.edu}

\address[Gatech]{School of Aerospace Engineering,
Georgia Institute of Technology, Atlanta, GA 30332 USA}  % Please supply                                              
\address[CMU]{Department of Mechanical Engineering, Carnegie Mellon University, Pittsburgh, PA 15213 USA}             % full addresses

\begin{keyword}                           
Stochastic control, covariance control, hybrid systems              
\end{keyword}                             

\begin{abstract}  
This work addresses the problem of optimally steering the state covariance of a linear stochastic system from an initial to a target, subject to hybrid transitions. The nonlinear and discontinuous jump dynamics complicate the control design for hybrid systems. Under uncertainties, stochastic jump timing and state variations further intensify this challenge. This work aims to regulate the hybrid system's state trajectory to stay close to a nominal deterministic one, despite uncertainties and noises. We address this problem by directly controlling state covariances around a mean trajectory, and this problem is termed the Hybrid Covariance Steering (H-CS) problem. The jump dynamics are approximated to the first order by leveraging the Saltation Matrix. When the jump dynamics are nonsingular, we derive an analytical closed-form solution to the H-CS problem. For general jump dynamics with possible singularity and changes in the state dimensions, we reformulate the problem into a convex optimization over path distributions by leveraging Schrodinger's Bridge duality to the smooth covariance control problem. The covariance propagation at hybrid events is enforced as equality constraints to handle singularity issues. The proposed convex framework scales linearly with the number of jump events, ensuring efficient, optimal solutions. This work thus provides a computationally efficient solution to the general H-CS problem. Numerical experiments are conducted to validate the proposed method.
 
\end{abstract}

\end{frontmatter}

\input{introduction}

\input{preliminaries}

\input{problem_formulation}

\input{main_results}

\input{experiments}

\input{conclusion}

\input{Appendix}

% \begin{ack}                               % Place acknowledgements
% Partially supported by the Roman Senate.  % here.
% \end{ack}

\bibliographystyle{plain}        % Include this if you use bibtex 
\bibliography{autosam}           % and a bib file to produce the 
                                 % bibliography (preferred). The
                                 % correct style is generated by
                                 % Elsevier at the time of printing.

%\begin{thebibliography}{99}     % Otherwise use the  
                                 % thebibliography environment.
                                 % Insert the full references here.
                                 % See a recent issue of Automatica 
                                 % for the style.
%  \bibitem[Heritage, 1992]{Heritage:92}
%     (1992) {\it The American Heritage. 
%     Dictionary of the American Language.}
%     Houghton Mifflin Company.
%  \bibitem[Able, 1956]{Abl:56}
%     B.~C.~Able (1956). Nucleic acid content of macroscope. 
%     {\it Nature 2}, 7--9. 
%  \bibitem[Able {\em et al.}, 1954]{AbTaRu:54}   
%     B.~C. Able, R.~A. Tagg, and M.~Rush (1954).
%     Enzyme-catalyzed cellular transanimations.
%     In A.~F.~Round, editor, 
%     {\it Advances in Enzymology Vol. 2} (125--247). 
%     New York, Academic Press.
%  \bibitem[R.~Keohane, 1958]{Keo:58}
%     R.~Keohane (1958).
%     {\it Power and Interdependence: 
%     World Politics in Transition.}
%     Boston, Little, Brown \& Co.
%  \bibitem[Powers, 1985]{Pow:85}
%     T.~Powers (1985).
%     Is there a way out?
%     {\it Harpers, June 1985}, 35--47.

%\end{thebibliography}

% \appendix
% \section{A summary of Latin grammar}    % Each appendix must have a short title.
% \section{Some Latin vocabulary}         % Sections and subsections are supported  
%                                         % in the appendices.
\end{document}

%% file: introduction.tex
\section{Introduction}
\label{sec:introduction}
Uncertainties are intrinsic in control systems, making feedback control crucial for minimizing their impact on system behavior. A natural objective is to maintain the system's performance close to the desired outcome despite these uncertainties. The covariance control paradigm \cite{HotSke87,IwaSke92,XuSke92} is a direct approach to managing uncertainties. The initial goal was to control a linear time-invariant stochastic system to achieve a specified steady-state covariance. More recently, covariance control theory was extended to the finite horizon control setting \cite{CheGeoPav15a,CheGeoPav15b,CheGeoPav15c,CheGeoPav17a,Bak16a,Bak18,yu2023covariance}, where the objective is to steer the state covariance of a linear dynamic system from an initial value to a desired target one. The covariance control framework has been effectively applied across a range of applications \cite{ZhuGriSke95,RidTsi18,OkaTsi19,okamoto2019optimal,Kal02,yu2023gaussian,yu2023stochastic}. However, the above paradigms and applications constrained their formulation within smooth dynamics. 

The co-existence of continuous state evolutions (also known as \textit{flows}) and discrete time events (also known as \textit{resets} or \textit{jump dynamics}) represents a wide range of real-world phenomena, such as power circuits \cite{giaouris2008stability} and bifurcations \cite{jiang2017grazing}. Examples of these \textit{Hybrid Dynamical Systems} in robotics include rigid body dynamics with contact \cite{sreenath2011compliant, posa2014direct}. The instantaneous and discontinuous jump dynamics through hybrid events make control, planning, and estimation challenging for hybrid systems. The limit cycles of hybrid systems exhibit small Regions of Attractions (ROA) and are difficult to stabilize around a reference `orbit' \cite{goswami1996limit, manchester2011transverse, manchester2011stable}. The optimal control of deterministic hybrid systems has been studied extensively \cite{goebel2009hybrid, heemels2012periodic}. Optimization-based formulations often treat hybrid transitions as complementary constraints in an optimization program \cite{posa2014direct, patel2019contact}. 

One way of addressing the complexity of jump dynamics is to approximate them as linear transitions. The most direct way is by taking the Jacobian with respect to the pre-event states and controls. However, this approximation is not accurate due to the instantaneous nature of the jumps. A more accurate linear approximation is the \textit{Saltation Matrix} (also known as the \textit{Jump Matrix}) \cite{filippov2013differential, leine2013dynamics, kong2023saltation}, accounting for the total variation at hybrid events caused by impact timing and jump dynamics. Covariance propagation and the Riccati equations at the hybrid events can also be approximated using the Saltation Matrix \cite{kong2023saltation}.

The linear approximation of jump dynamics allows for the application of smooth linear optimal control methods to hybrid systems. The Linear Quadratic Regulator (LQR) or Linear Quadratic Gaussian (LQG), is one classical paradigm \cite{kailath1980linear}. LQR has previously been extended to periodic systems in \cite{carnevale2014linear, possieri2020linear}. However, only linear time-invariant (LTI) systems with time-triggered switching hybrid dynamics are considered. Based on the Saltation Matrix approximation, iterative methods for nonlinear smooth problems were extended to hybrid systems in \cite{kong2021ilqr, kong2023hybrid}, where the event-based hybrid iterative-LQR (H-iLQR) was introduced as the counterpart to the smooth Differential Dynamic Programming (DDP) \cite{theodorou2010stochastic, tassa2014control}. 

Smooth covariance steering differs from the LQG control design in that the uncertainties represented by the covariances for linear systems are guaranteed to be steered to a target value instead of relying on a terminal cost function design in LQG. The formulations and the solution techniques are thus significantly different. In LQG, the problem resorts to solving the Riccati equations starting from a terminal quadratic loss, while the covariance control problem solves a coupled Lyapunov and Riccati equation pair such that the terminal constraint on the state covariance for the closed-loop system is satisfied.

This work designs controllers that regulate the state covariances of linear stochastic systems around a mean trajectory, subject to hybrid transitions. Our main contribution is showing that this problem has a closed-form solution or can be solved efficiently via a small-size convex optimization, within the covariance steering paradigm. We approximate the jump dynamics to the first order using the Saltation Matrix. The problem has a closed analytical solution similar to the smooth covariance control settings for invertible jumps. For general jump dynamics, we formulate the problem into a convex optimization, which allows for singular covariance at post-impact instances \cite{ciccone2020regularized}. This work extends the smooth covariance steering paradigm \cite{CheGeoPav15a,CheGeoPav15b,CheGeoPav15c,CheGeoPav17a} to hybrid dynamics settings, which we believe has a vast range of potential applications in control engineering and robotics \cite{iqbal2020provably, gu2018exponential}. 

The motivation for this work also stems from the growing complexity of control tasks in robotics and autonomous systems, with hybrid dynamics—such as legged locomotion and robotic manipulation with intermittent contact. These tasks require efficient and principled methods to control the uncertainties. In these scenarios, managing uncertainty is critical to ensuring reliable system performance, particularly in stochastic environments with hybrid state transitions. Existing methods often encounter challenges in scalability and computational efficiency, especially when dealing with contact-rich systems. By providing a closed-form solution or a small size convex optimization formulation for steering state covariances in hybrid systems, this work aims to offer a more robust and computationally efficient approach for such applications. 

The works that are most related to this work are the smooth covariance steering theory \cite{CheGeoPav15a, CheGeoPav14e, CheGeoPav17a}, and the hybrid iLQR control designs \cite{possieri2020linear, kong2021ilqr, kong2023hybrid}. Previous efforts to manage uncertainties in contact-rich systems include optimization-based methods introducing stochastic complementary constraints to characterize uncertain contacts \cite{tassa2011stochastic, drnach2021robust}. However, these methods lift contact constraints into the cost functions and are not always guaranteed to be met. In \cite{zhu2024convergent}, the authors added a penalty term in the terminal covariance within the H-iLQR framework to control the uncertainty levels for hybrid systems. 

The rest of the paper is organized as follows: Section \ref{sec:preliminaries} introduces the smooth covariance steering problem. Section \ref{sec:problem_formulation} presents the H-CS problem formulation. Section \ref{sec:E_invertible} presents a closed-form solution for the invertible linear reset maps, and Section \ref{sec:E_rectangular} presents a convex optimization paradigm for general linear reset maps. Section \ref{sec:experiments} contains numerical examples, followed by a conclusion in Section \ref{sec:conclusion}.

%% file: preliminaries.tex
\section{Preliminaries}
\label{sec:preliminaries}
\subsection{Smooth Linear Covariance Steering}
Smooth covariance steering optimally steers an initial Gaussian distribution to a target one at the terminal time while minimizing the control energy and a quadratic state cost. We denote $X(t) \in \mathbb{R}^{n}$ as the continuous-time state, and $u(t) \in \mathbb{R}^{m}$ as the control input. The smooth covariance steering problem is defined as the following optimal control problem \cite{CheGeoPav15a, CheGeoPav15b, CheGeoPav17a}
\begin{subequations}\label{eq:covcontrollinear}
\begin{align}
    \!\!\!\!\!\min_{u(t)} \; \mathcal{J} & =  \mE \left\{\int_{0}^{T} [\|u(t)\|^2 + X'(t) Q(t) X(t)]dt\right\}
    \label{eq:cov_steering_obj}
    \\ 
    \!\!\!\!dX(t) & = A(t) X(t) dt + B(t) (u(t) dt+ \sqrt{\epsilon} d W(t)) 
    \label{eq:cov_constraint}
    \\ 
    \!\!\!\!X(0) & \sim \cN(m_0, \Sigma_0),\quad X(T) \sim \cN(m_T, \Sigma_T), \label{eq:cov_constraint_j}
\end{align}
\end{subequations}
where $W(t)$ is a standard Wiener process with intensity $\epsilon$, $Q(t)$ is a running cost, and $\mathcal{N}(m_0, \Sigma_0)$ and $\mathcal{N}(m_T, \Sigma_T)$ represent the designated initial and terminal time Gaussian distributions. Thanks to linearity, the mean and covariance can be controlled separately.

Assume $\Pi(t)$ satisfies the Riccati equations
\begin{align}
    \label{eq:riccati_smooth}
    &-\dot\Pi(t) = \\
    &A'(t)\Pi(t) + \Pi(t) A(t)  -  \Pi(t) B(t) B'(t) \Pi(t) + Q(t),\nonumber
\end{align}
then the optimal control to \eqref{eq:covcontrollinear} has a feedback form 
    \begin{equation}
    \label{eq:optimal_feedback}
        u(t)^* = - B'(t)\Pi(t) X(t), \; j=1,2,
    \end{equation}
and the evolution of the controlled covariance function $\Sigma(t)$ is governed by the Lyapunov equations
\begin{equation}
\label{eq:covariance_evolution_smooth}
\begin{split}
    \!\!\!\!\dot{\Sigma}(t) = ( & A(t)-B(t)B'(t)\Pi(t))\Sigma(t) + \Sigma(t)(A(t)-
    \\
    \!\!\!\!\!&B(t)B'(t)\Pi(t))' + B(t)B'(t).
\end{split}
\end{equation}
Define an intermediate variable, $H(t)$ as
\begin{equation}
\label{eq:defn_H}
    H(t) \triangleq \Sigma^{-1}(t) - \Pi(t).
\end{equation}
Consider \eqref{eq:riccati_smooth}, \eqref{eq:covariance_evolution_smooth}, a direct calculation gives
\begin{align}
\label{eq:riccati_H}
    &-\dot{H}(t) = \\
    &A'(t)H(t) + H(t)A(t) + H(t)B(t)B'(t)H(t) - Q(t).\nonumber
\end{align}
The optimal solution to \eqref{eq:covcontrollinear} is given by the solution to the two Riccati equations \eqref{eq:riccati_smooth}, \eqref{eq:riccati_H} coupled by boundary conditions 
\begin{subequations}
\label{eq:LQschrodinger}
    \begin{align}
    \label{eq:LQschrodinger3}
    \epsilon \Sigma_0^{-1}&=\Pi(0) + H(0)
    \\
    \label{eq:LQschrodinger4}
    \epsilon \Sigma_T^{-1}&= \Pi(T) + H(T).
    \end{align}
    \end{subequations}
Denote the Hamiltonian systems for $(A, B)$ as 
\begin{equation}
\label{eq:defn_M}
   M(t)= \left[
    \begin{matrix}A(t) & -B(t) B'(t)\\-Q(t)& -A'(t)\end{matrix}
    \right] \in \mR^{2n \times 2n}
\end{equation}
with the transition kernels 
\begin{equation}
\label{eq:transition_kernel_M}
\Phi(t,s) \triangleq \begin{bmatrix}
    \Phi_{11}(t,s) & \Phi_{12}(t,s)
    \\
    \Phi_{21}(t,s) & \Phi_{22}(t,s)
\end{bmatrix},
\end{equation}
satisfying $\partial \Phi(t,s)/\partial t = M(t) \Phi(t,s)$ with $\Phi(s,s)=I$. For all the transition kernels that appear hereafter, omitting the time variables $(t,s)$ represents the transition from $t_0$ to $T$. e.g., $\Phi_{11} \triangleq \Phi(T,0)_{11}$. We use $\Psi$ to represent the reverse transition kernel, i.e., $\Psi \triangleq \Phi(0,T) = \Phi(T,0)^{-1}$. $\Psi$ always exists since $\Phi$ is never singular. $\Psi_{ij}$ denotes the blocks in the reversed transition kernel $\Psi$ similar to \eqref{eq:transition_kernel_M}.

The Riccati equations \eqref{eq:riccati_smooth}, \eqref{eq:riccati_H} with boundary conditions \eqref{eq:LQschrodinger} have a closed-form solution that is based on the properties of $\Phi(t,s)$, which we introduce in the following Lemma \ref{lem:state_transition}. 
\begin{lemma}
\label{lem:state_transition}
    (\cite{CheGeoPav17a}, Lemma 3 ). The entries of the state transition kernel $\Phi(t,s)$ satisfy 
    \begin{subequations}
    \label{eq:lemma_statetransition}
        \begin{align}
            \Phi_{11}(t,s)'\Phi_{22}(t,s) - \Phi_{21}(t,s)'\Phi_{12}(t,s) &= I,
            \label{eq:lemma_statetransition_1}
            \\
            \Phi_{12}(t,s)'\Phi_{22}(t,s) - \Phi_{22}(t,s)'\Phi_{12}(t,s) &= 0,
            \label{eq:lemma_statetransition_2}
            \\
            \Phi_{21}(t,s)'\Phi_{11}(t,s) - \Phi_{11}(t,s)'\Phi_{21}(t,s) &= 0,
            \label{eq:lemma_statetransition_3}
            \\
            \Phi_{11}(t,s)\Phi_{22}(t,s)' - \Phi_{12}(t,s)\Phi_{21}(t,s)' &= I,
            \label{eq:lemma_statetransition_4}
            \\
            \Phi_{12}(t,s)\Phi_{11}(t,s)' - \Phi_{11}(t,s)\Phi_{12}(t,s)' &= 0,
            \label{eq:lemma_statetransition_5}
            \\
            \Phi_{21}(t,s)\Phi_{22}(t,s)' - \Phi_{22}(t,s)\Phi_{21}(t,s)' &= 0,
            \label{eq:lemma_statetransition_6}
        \end{align}
    \end{subequations}
    and the entries $\Phi_{12}(t,s)$ and $\Phi_{11}(t,s)$ are invertible for all $t$, and $(\Phi_{12}(t,t_0)^{-1}\Phi_{11}(t,t_0))^{-1}$ is monotonically decreasing in the positive definite sense with left limit $0$ as $t \searrow t_0$.
\end{lemma}
We summarize the solution to the smooth covariance steering problem in the following proposition. 
\begin{proposition}
\label{thm:smooth_result}
The Riccati equations \eqref{eq:riccati_smooth} and \eqref{eq:riccati_H} with coupled boundary conditions \eqref{eq:LQschrodinger} have a unique solution, determined by the
initial values
\begin{subequations}\label{eq:initial}
\begin{align}
    \Pi(0) =& \;\frac{\epsilon\Sigma_0^{-1}}{2}- \Phi_{12}^{-1}\Phi_{11} - 
    \label{eq:initial_a}
    \\
    \nonumber & \; \Sigma_0^{-\frac{1}{2}}\!\left(\frac{\epsilon^2I}{4}+\Sigma_0^{\frac{1}{2}}\Phi_{12}^{-1}\Sigma_T (\Phi_{12}')^{-1}\Sigma_0^{\frac{1}{2}}\right)^{\frac{1}{2}}\!\Sigma_0^{-\frac{1}{2}},
    \\
    H(0)  =& \; \epsilon\Sigma_0^{-1}-\Pi(0).
    \end{align}
\end{subequations}
\end{proposition}

\begin{remark}
\label{thm:remark_terminal_Pi}
The above conditions for $\Pi$ in Eq. \eqref{eq:initial_a} can also be given in the form of terminal values of $\Pi(T)$ as
    \begin{align}
    \label{eq:terminal_a}
    \begin{split}
        \!\!\!\! \Pi(T) = &\; \frac{\epsilon\Sigma_T^{-1}}{2}- \Psi_{12}^{-1}\Psi_{11} + 
        \\
        & \; \Sigma_T^{-\frac{1}{2}}\!\left(\frac{\epsilon^2I}{4}+\Sigma_T^{\frac{1}{2}}\Psi_{12}^{-1}\Sigma_0 (\Psi_{12}')^{-1}\Sigma_T^{\frac{1}{2}}\right)^{\frac{1}{2}}\!\Sigma_T^{-\frac{1}{2}},
    \end{split}
    \end{align}
    The proof is similar to that in \cite{CheGeoPav17a}, Theorem 1, and is omitted.  
\end{remark}

\subsection{Schr\"odinger Bridge Problem}
Covariance steering is closely related to Schrodinger's Bridge Problem (SBP) of finding the closest distribution to a prior one that connects two marginals \cite{CheGeoPav14e, CheGeoPav15a}. We use $\mP$ to denote the measure induced by the controlled process \eqref{eq:cov_constraint} with marginal covariance constraints \eqref{eq:cov_constraint_j}. Letting $u(t) \equiv 0$, we define the prior distribution as the one induced by the uncontrolled process
\[
dX(t) = A(t) X(t) dt + \sqrt{\epsilon}B(t) dW(t),
\]
whose induced path measure is denoted as $\mQ$. Denote the state-related running cost function as
\begin{equation*}
    \mathcal{L} \triangleq \int_{0}^{T} X'(t)Q(t)X(t) dt,
\end{equation*} 
by Girsanov's theorem \cite{Gir60}, the expected control energy equals to the relative entropy between $\mP$ and $\mQ$, i.e.,
\begin{equation}
\label{eq:Girsanov}
    \mE \left\{\int_{0}^{T} \frac{1}{2\epsilon} \|u(t)\|^2 dt \right\} = \mE_{\mP} \left[ \log\left( \frac{d\mP}{d\mQ} \right) \right],
\end{equation}
leading to the following equivalent objective function to \eqref{eq:covcontrollinear} 
\begin{subequations}
\label{eq:obj_KL}
\begin{align}
     \mathcal{J} &= \mE_{\mP} \left[ 2 \epsilon \log\left( \frac{d\mP}{d\mQ} \right) + \mathcal{L} \right] 
    \label{eq:obj_measure_1}
    \\
     &\propto \mE_{\mP}  \left[ \epsilon \log d\mP - \epsilon \log \left( \frac{d\mQ e^{-\mathcal{L} / 2\epsilon}}{\mE_{\mQ} \left[ e^{-\mathcal{L}/ 2 \epsilon} \right] } \right) \right]  \nonumber
    \\
    &= {\rm KL} \left( \mP \parallel \mP^* \right), 
    \label{eq:obj_measure_2}
\end{align}
\end{subequations}
where in the second equation we add and then subtract a fixed term $\epsilon \log \mE_{\mQ} [e^{-\frac{\mathcal{L}}{2 \epsilon}}]$, and the measure $\mP^*$ is defined as
\begin{equation}
\label{eq:defn_dmP_star}
    d\mP^* \triangleq \frac{e^{-\mathcal{L}/ 2 \epsilon}}{\mE_{\mQ} \left[ e^{-\mathcal{L} / 2 \epsilon} \right]} d\mQ.
\end{equation}
The smooth covariance steering problem \eqref{eq:covcontrollinear} is therefore equivalent to a KL-divergence minimization problem.

%% file: problem_formulation.tex
\section{Problem Formulation}
\label{sec:problem_formulation}
This section presents the main problem formulation.
\subsection{Hybrid Systems and Saltation Approximation}
A hybrid dynamical system contains different modes with different smooth flow and discrete jump dynamics between them. It is usually defined by a tuple \cite{grossman1993hybrid, johnson2016hybrid} $\mathcal{H} \coloneqq \{ \mathcal{I}, \mathcal{D}, \mathcal{F}, \mathcal{G}, \mathcal{R} \}.$
The set 
\begin{equation}
    \mathcal{I} \coloneqq \{ I_1, I_2, \dots, I_{N_I}\} \subset \mathbb{N}
\label{eq:hybrid_definitions}
\end{equation}
is a finite set of \textit{modes}, $\mathcal{D}$ is the set of continuous \textit{domains} containing the state spaces, with domain $D_{j}$ for mode $I_j$, $\mathcal{F}$ is the set of \textit{flows}, consisting of individual flow $F_{j}$ that describes the smooth dynamics in mode $I_j$. $\mathcal{G} \subseteq \mathcal{D}$ denotes the set of \textit{guards} triggering the resets. 

We denote $X_j(t) \in \mathbb{R}^{n_j}$ as the continuous-time state in mode $I_j$, and $u_j(t) \in \mathbb{R}^{m_j}$ as the state-dependent control input. In mode $I_j$, the flow $F_j$ is considered to be a linear stochastic process
\begin{equation}
\label{eq:dyn_smooth_j}
    d X_j(t) = (A_{j}(t)X_j(t)  + B_{j}(t) u_j(t)) d t + \sqrt{\epsilon}B_j(t) d W_j(t)), 
\end{equation}
where $A_j(t) \in \mR^{n_j \times n_j}, B_j(t) \in \mR^{n_j \times m_j}$ are the system matrices in mode $I_j$, and $dW_j(t) \in \mR^{m_j}$ denotes a standard Wiener process with noise intensity $\epsilon$. All the above variables are \textit{mode-dependent}. A transition from mode $I_j$ to mode $I_k$ happens at state $X_j(t^-)$ and time $t^{-}$ if the guard condition is satisfied, i.e.,
\begin{equation*}
    g_{jk}(t^{-}, X_j(t^-)) \leq 0, \; X_j(t^-) \in I_j.
\end{equation*}
The guard function triggers the time instances and states on which the jump dynamics events happen. At the event $t^-$, an instant jumping dynamics is applied to the system. From mode $I_j$ to mode $I_k$, the reset map is defined as
\begin{equation*}
    t^+ = t^-, \;\; X_k(t^+) = R_{jk} (X_j(t^-)).
\end{equation*} 
The reset map $R_{jk}$ is nonlinear and instantaneous in general. Saltation Matrices \cite{kong2023saltation} provide precise linear approximations for the nonlinear reset functions. The Saltation Matrix $\Xi_{jk} \in \mR^{n_k \times n_j}$ from mode $I_j$ to $I_k$ is defined as 
\begin{equation}
\label{eq:saltation}
    \Xi_{jk} \triangleq \partial_x R_{jk} + \frac{(F_{k} - \partial_x R_{jk} \cdot F_{j} - \partial_t R_{jk})\partial_x g_{jk}}{\partial_t g_{jk} + \partial_x g_{jk} \cdot F_{j}},
\end{equation}
where $\partial_x(\cdot)$ and $\partial_t(\cdot)$ denote the partial derivatives in state and time, respectively. The dynamics of a perturbed state $\delta X_t$ at hybrid transitions can be approximated as
\begin{equation}
\label{eq:saltation_approximation_dyn}
    \delta X(t^+) \approx \Xi_{jk} \delta X(t^-).
\end{equation}
Built on this approximation, controller design tasks based on iterative linearization of the dynamics have been applied to robot locomotion and estimation tasks for hybrid systems \cite{kong2023hybrid, kong2021salted}. Motivated by this linear approximation, this work considers \textit{saltation} dynamics for the stochastic state variable $X(t)$ around its mean at the hybrid events 
\begin{equation}
\label{eq:hybrid_dyn}
    t^+ = t^-, \;\; X_k(t^+) = \Xi_{jk} X_j(t^-),
\end{equation} 
where $n_j, n_k$ are the dimensions of the state space before and after the jump. With the above definitions of saltation transition and linear stochastic flow, the hybrid dynamics considered in this work are well-defined. 

\subsection{Covariance Steering with Hybrid Transitions}
{\em Assumptions and Notations.}
We consider our optimal control problems in the time window $[0, T]$. Without loss of generality, we only consider \textit{one} jump event at $t=t^-$, which separates the total time window into $2$ pieces. Denote the initial and terminal time in the two windows as $t_0^j$ and $t_f^j$ respectively, with $t_0^1 = 0, t_f^1 = t^-$ and $t_0^2=t^+, t_f^2=T$.
Denote the two modes before and after the jump event as $I_1$ and $I_2$, where the pre-event states $X_1(t) \in \mathbb{R}^{n_1}, \forall t \leq t^-$ and the post-event states $X_2(t) \in \mathbb{R}^{n_2}, \forall t \geq t^+$. We use 
\[
\Xi \triangleq \Xi_{12} \in \mR^{n_2\times n_1}
\]
to represent the saltation transition from mode $I_1$ to $I_2$ at $t^-=t^+$.

The covariance steering for systems with saltation transitions is formulated as follows 
\begin{subequations}\label{eq:formulation_main}
\begin{align}
    & \hspace*{-1.4cm} \min_{u_j(t)} \mathcal{J}_H  \triangleq \mE  \left\{ \int_0^{T} [ \|u_j(t)\|^2 + X_j'(t) Q_j(t) X_j(t)]dt \right\} \nonumber
    \\ 
    \!dX_1  &= A_1(t) X_1 dt +  B_1(t)( u_1 dt + \! \sqrt{\epsilon} dW_1),
    \label{eq:hybrid_lin_sys1} 
    \\ 
    \!X_2(t^+)  &=  \Xi X_1(t^-), 
    \label{eq:hybrid_lin_sys2} 
    \\ 
    \!dX_2 &=  A_2(t) X_2 d t  +   B_2(t)( u_2 dt + \! \sqrt{\epsilon}dW_2),
    \label{eq:hybrid_lin_sys3} 
    \\   
    \!X_1(0) &\sim  \cN(m_0, \Sigma_0), \quad \!\! X_2(T) \sim \cN(m_T, \Sigma_T),\label{eq:formulation_covariance_constraints} 
\end{align}
\end{subequations}
where $m_0 \in \mR^{n_1}, \Sigma_0 \in \mR^{n_1\times n_1}$ are the designated initial state mean and covariance in mode $I_1$, and $m_T\in \mR^{n_2}, \Sigma_T \in \mR^{n_2\times n_2}$ are the terminal state mean and a designated target state covariance in mode $I_2$, respectively. 

The other quantities and matrices that will appear in this section, i.e., ($X(t),u(t),A(t),B(t),Q(t),W(t), \Sigma(t), \Pi(t)$), all have the same definitions as in the smooth problem \eqref{eq:covcontrollinear} with a mode dependent notation $j\in \{1,2\}$ to indicate the mode $I_j$ they are in. The jump event time $t^-$ is determined by the controlled mean trajectory, decomposing the objective function $\mathcal{J}_H(u)$ into $3$ parts
\begin{equation}
\label{eq:J_decompose}
\mathcal{J}_H(u) \coloneqq \mathcal{J}_1(u) + \mathcal{J}_2(u) + \mathcal{J}_3(u),
\end{equation} 
where $\mathcal{J}_{j\in\{1,2\}}$ are the costs for the $2$ smooth problems 
\begin{subequations}
% \label{eq:defn_J1_J2}
    \begin{align}
      & \! \mathcal{J}_{j}(u) = \mE \left\{\int_{t_0^j}^{t_f^j} [\|u_j(t)\|^2 + X_j'(t) Q_j(t) X_j(t)]dt \; + \right. \nonumber
        \\
        & \left. \quad \vphantom{}X_j'(t_f^j)\Pi_j(t_f^j) X_j(t_f^j) - X_j'(t_0^j)\Pi_j(t_0^j)X_j(t_0^j) \right\} 
        \\
        &\! = \mE \left \{\int_{t_0^j}^{t_f^j} [\|u_j(t)\|^2 + X_j'(t) Q_j(t) X_j(t)]dt + \right. \nonumber
        \\
        & \quad\quad\quad\quad d\left(X_j'(t)\Pi_j(t) X_j(t)\right) \Big \},
    \end{align}
\end{subequations} 
and $\mathcal{J}_3$ is the cost at the jumping event, which is
\begin{subequations}
\begin{align*}
    \mathcal{J}_3 &= \mE \left\{ \int_{t^-}^{t^+} d \left( X_j(t)'\Pi_j(t)X_j(t)\right)\right\} 
    \\
   & \!\! = \mE \left[X_2'(t^+)\Pi_2(t^+)X_2(t^+) - X_1'(t^-)\Pi_1(t^-)X_1(t^-) \right].
\end{align*}
\label{eq:defn_J3}
\end{subequations}
Assume that $\Pi_1(t^-)$ and $\Pi_2(t^+)$ satisfy the discrete-time Riccati equation
\begin{equation}
    \label{eq:riccati_hybrid}
    \Xi'\Pi_2(t^+)\Xi = \Pi_1(t^-),
\end{equation} 
then we have 
\begin{equation*}
% \label{eq:J3_expanded}
    \mathcal{J}_3 \equiv 0.
\end{equation*}
Consider the Riccati equations \eqref{eq:riccati_smooth}, and expand the term
\[
d(X_j'(t)\Pi_j(t) X_j(t))
\]
by Ito's rule in $J_{j\in \{1,2\}}$, we have 
\begin{align}
\label{eq:J_expanded}
\mathcal{J}_H(u) = \mE & \left\{ \sum_{j} \int_{t_0^j}^{t_f^j} [\|u_j(t) + B_j'(t)\Pi_j(t) X_j(t) \|^2 + \right.  \nonumber
\\
& \Pi_j(t) B_j(t) B_j'(t)\Pi_j(t) ]dt \Biggr\}.
\end{align}
Minimizing the quadratic function \eqref{eq:J_expanded} for $u_j(t)$, the optimal control for the overall problem is given in the feedback form as in \eqref{eq:optimal_feedback}, and the controlled covariance evolution follows \eqref{eq:covariance_evolution_smooth}.
At the jump event, the covariance evolves as
\begin{equation}
\label{eq:covariance_evolution_hybrid}
\begin{split}
    \!\!\!\!\!\Sigma_2(t^+) &= \mE \left[ (X_2(t^+) - x_2(t^+)) (X_2(t^+) - x_2(t^+))' \right] 
    \\
    &= \Xi \Sigma_1(t^-) \Xi'.
\end{split}
\end{equation}
We re-state the formulation \eqref{eq:formulation_main} as the following Problem \ref{prob:overall_problem}.
\begin{problem}
    \label{prob:overall_problem}
    Solving problem \eqref{eq:formulation_main} is equivalent to solving the pair $(\Pi_j(t), \Sigma_j(t))$ satisfying \eqref{eq:riccati_smooth}, \eqref{eq:riccati_hybrid}, \eqref{eq:covariance_evolution_smooth}, \eqref{eq:covariance_evolution_hybrid}, with boundary conditions 
    \begin{equation}
\label{eq:covariance_constraints}
    \Sigma_1(0) = \Sigma_0, \; \Sigma_2(T) = \Sigma_T.
\end{equation}
The optimal solution is in the feedback form \eqref{eq:optimal_feedback}. 
\end{problem}

From the decomposition \eqref{eq:J_expanded}, we see that the H-CS problem can be decomposed into $2$ smooth covariance steering problems, in $[0, t^-]$ and $[t^+, T]$, with marginal covariance at $t=0, t=t^-, t=t^+, t=T$. The propagation of $\Sigma$ and $\Pi$ at $(t^-, t^+)$ satisfy \eqref{eq:covariance_evolution_hybrid} and \eqref{eq:riccati_hybrid}. Consider a feedback policy \eqref{eq:optimal_feedback} in $[0, t^-]$ that optimally steers the state covariance from $\Sigma_0$ to $\Sigma_1(t^-) = \Sigma^-$, then the covariance $\Sigma_2(t^+) = \hat\Sigma^{+}=\Xi\Sigma^- \Xi'$ serves as the initial condition for the smooth covariance steering problem in $[t^+, T]$. For any $\Sigma^-$, we can find the optimal controller that steers $\Sigma^+$ to $\Sigma_T$ and minimizes $J_2(u)$. The optimal solution to the overall H-CS requires that $\Pi_1(t^-)$ and $\Pi_2(t^+)$ satisfy \eqref{eq:riccati_hybrid}.

The $\Pi_2(t^+)$ is obtained by replacing $\Sigma_0$ by $\Xi \Sigma^- \Xi'$ in \eqref{eq:initial_a}. Similarly, the optimal $\Pi_1(t^-)$ is a function of $\Sigma^-$ by replacing $\Sigma_T$ with $\Sigma^-$ in \eqref{eq:initial_a}. The covariances $\Sigma^-$ and $\Sigma^+$ thus become variables in the overall hybrid covariance control problem. By equating $\Pi_2(t^+)$ with $\Xi'\Pi_1(t^-)\Xi$, we establish the condition for the optimal $\Sigma^-$, $\Sigma^+$, $\Pi_1(t^-)$ and $\Pi_2(t^+)$ for the overall problem. It turns out that the shape and rank properties of $\Xi$ pose technical issues in solving the problem in general. In the sequel, in Section \ref{sec:E_invertible}, we first solve the problem in analytical form for an invertible $\Xi$. In Section \ref{sec:E_rectangular}, for general $\Xi$ matrices, we solve the problem numerically by formulating it into a convex optimization.

%% file: main_results.tex
\section{Results for Invertible Jump Dynamics}
\label{sec:E_invertible} 
This section derives the solution to Problem \ref{prob:overall_problem} for linear stochastic systems with an \textit{invertible} hybrid transition $\Xi$. The smooth evolutions for the matrices $\Sigma_j(t), H_j(t)$, and $\Pi_j(t)$ stay the same. At the hybrid time, we have 
\begin{subequations}
\label{eq:hybrid_transition_H}
\begin{align}
     H_2(t^+) &= \Sigma_2^{-1}(t^+) - \Pi_2(t^+)
     \nonumber
    \\
    &= (\Xi\Sigma_1(t^-)\Xi')^{-1} - \Pi_2(t^+)
    \\
     H_1(t^-) &= \Sigma_1^{-1}(t^-) - \Pi_1(t^-)
     \nonumber
     \\
     &= \Sigma_1^{-1}(t^-) - \Xi'\Pi_2(t^+)\Xi.
\end{align}
\end{subequations}
The jumps of $\Pi_j(t)$ and $H_j(t)$ at $t=t^-$ are thus
\begin{subequations}
\begin{align}
    \Pi_2(t^+) &= (\Xi')^{-1} \Pi_1(t^-) \Xi^{-1}, \label{eq:riccati_Pi_invertibleE}
    \\
    H_2(t^+) &= (\Xi')^{-1} H_1(t^-) \Xi^{-1}.\label{eq:riccati_H_invertibleE}
\end{align}
\end{subequations}
With $H_j(t)$, the boundary conditions \eqref{eq:covariance_constraints} is equivalently
\begin{subequations}
\label{eq:coupled_boundary_Pi_H}
\begin{align}
    \Sigma^{-1}_1(0) &= \Pi_1(0) + H_1(0), 
    \\
    \Sigma^{-1}_2(T) &= \Pi_2(T) + H_2(T).
\end{align}
\end{subequations}
We summarize the above into the following Problem \ref{prob:prob_invertible_E}.
\begin{problem}
    \label{prob:prob_invertible_E}
    For invertible jump dynamics $\Xi$, solving Problem \ref{prob:overall_problem} is equivalently solving the pair $(\Pi_j(t), H_j(t))$ for equations \eqref{eq:riccati_smooth}, \eqref{eq:riccati_hybrid}, \eqref{eq:riccati_H}, \eqref{eq:hybrid_transition_H}, with boundary conditions \eqref{eq:coupled_boundary_Pi_H}. 
\end{problem}

We first define the transition kernel $\Phi^\Xi$ for the jump dynamics as
\begin{equation}
\label{eq:defn_Phi_minus_invertible}
    \Phi^{\Xi}(t^+, t^-) = 
    \begin{bmatrix}
     \Phi^{\Xi}_{11} & \Phi^{\Xi}_{12}
     \\
     \Phi^{\Xi}_{21} & \Phi^{\Xi}_{22}
    \end{bmatrix}
    \triangleq 
    \begin{bmatrix}
     \Xi & 0
     \\
     0 & (\Xi')^{-1}
    \end{bmatrix}.
\end{equation}
For time $t>t^-$, we define the hybrid transition kernel $\Phi^H(t,0)$ as
\begin{equation*}
% \label{eq:defn_PhiH}
    \Phi^H(t,0) \triangleq 
    \begin{cases}
    \Phi_1(t, 0) , & \text{if}\ t \leq t^- 
      \\
    \Phi_2(t, t^+) \Phi^{\Xi}(t^+,t^-) \Phi_1(t^-, 0), & \text{if}\ t > t^-
    \end{cases}
\end{equation*}
where $\Phi_j(t,s)$ are defined in \eqref{eq:transition_kernel_M}. 

\begin{lemma}
\label{lem:lemma_hybrid_Phi_invertible_E}
Denote $[X_H'(t), Y_H'(t)]' = \Phi^H(t,0) [X_0', Y_0']'$, then 
\begin{equation}
\label{eq:sol_riccati_XY}
    \Pi^H_t \triangleq Y_H(t) X^{-1}_H(t)
\end{equation}
is a solution to the equations \eqref{eq:riccati_smooth} and \eqref{eq:riccati_Pi_invertibleE} for all $t \in [0,T]$. Denote $[\hat{X}_H'(t), \hat{Y}_H'(t)]' = \Phi^H(t,0) [\hat{X}_0', \hat{Y}_0']'$, the term
\begin{equation*}
% \label{eq:sol_riccati_H}
    H^H_t \triangleq -\hat{X}_H'(t)^{-1}\hat{Y}_H(t)
\end{equation*}
is a solution to equations \eqref{eq:riccati_H} and \eqref{eq:riccati_H_invertibleE} for all $t \in [0,T]$.

\begin{proof}
See Appendix \ref{sec:appendix_proof_Lemma_hybrid_Phi}.
\end{proof}
\end{lemma}
The following Lemma \ref{lem:state_transition_product} and Lemma \ref{lem:statetransition_hybrid_invertible_E} show that the hybrid transition kernel $\Phi^H$ preserves the properties in Lemma \ref{lem:state_transition}. 
\begin{lemma}
\label{lem:state_transition_product}
    If two state transitions, $\bar{\Phi}(t,s)$ and $\Tilde{\Phi}(t,s)$ both satisfy the equations \eqref{eq:lemma_statetransition_1} - \eqref{eq:lemma_statetransition_6}, then the product $\Phi(t,s) = \bar{\Phi}(t,s)\Tilde{\Phi}(t,s)$ also satisfies \eqref{eq:lemma_statetransition_1} - \eqref{eq:lemma_statetransition_6}.
    
    \begin{proof}
        See Appendix \ref{sec:appendix_lemma_product}.
    \end{proof}
\end{lemma}
\begin{lemma}
\label{lem:statetransition_hybrid_invertible_E}
The transition kernel $\Phi^H(t,0)$ satisfies the equations \eqref{eq:lemma_statetransition_1} - \eqref{eq:lemma_statetransition_6}, the entries $\Phi^H_{12}(t,s)$ and $\Phi^H_{11}(t,s)$ are invertible for all $t$, and $(\Phi^H_{11}(t,0))^{-1}\Phi^H_{12}(t,0)$ is monotonically decreasing in the positive definite sense with left limit $0$ as $t \searrow 0$.

\begin{proof}
    See Appendix \ref{sec:proof_statetransition_hybrid_invertible_E}. 
\end{proof}
\end{lemma}
Lemma \ref{lem:statetransition_hybrid_invertible_E} shows that $\Phi^H$ preserves the conditions in Lemma \ref{lem:state_transition}. We present the main results in Theorem \ref{thm:main}. 
\begin{theorem}[Main results for invertible jump dynamics]
\label{thm:main}
Denote the hybrid transition matrix as
\begin{align*}
    \begin{bmatrix}
    \Phi^H_{11} & \Phi^H_{12}
    \\
    \Phi^H_{21} & \Phi^H_{22}
    \end{bmatrix}
    \triangleq 
    \begin{bmatrix}
    \Phi^H_{11}(T,0) & \Phi^H_{12}(T,0)
    \\
    \Phi^H_{21}(T,0) & \Phi^H_{22}(T,0)
    \end{bmatrix}.
\end{align*}
The initial conditions for $(\Pi_j(t), H_j(t))$ that solves \eqref{eq:formulation_main} are 
\begin{subequations}
\label{eq:main_result_Pi0_H0}
    \begin{align}
         \Pi_1(0) &= \frac{1}{2}\Sigma_0^{-1} -(\Phi^H_{12})^{-1}\Phi^H_{11} - \nonumber
        \\
        \!\!\!\!\!\! & \!\!\!\!\!\!\!\!\! \Sigma_0^{-\frac{1}{2}} (\frac{I}{4} + \Sigma_0^{\frac{1}{2}} (\Phi^H_{12})^{-1} \Sigma_T ((\Phi^H_{12})')^{-1} \Sigma_0^{\frac{1}{2}})^{\frac{1}{2}} \Sigma_0^{-\frac{1}{2}},
        \\
        H_1(0) &= \Sigma_0 - \Pi_0.
    \end{align}
\end{subequations}
\begin{proof}
    The proof is similar to the proof of Theorem 1 in \cite{CheGeoPav17a}, based on Lemma \ref{lem:statetransition_hybrid_invertible_E}. Consider boundary conditions \eqref{eq:coupled_boundary_Pi_H}. Let $X_0=\hat{X}_0=I$ without affecting the value of $\Pi_1(0)$, we have
\begin{subequations}
    \begin{align}
        \Sigma_0^{-1} =&\; Y_0 - \hat Y_0',
        \label{eq:boundary_Y0_hatY0}
        \\
        \Sigma_T^{-1} =&\; Y_H(T) X_H(T)^{-1} - (\hat X_H(T)')^{-1}\hat Y_H(T)' \nonumber
        \\ 
        =&\; (\Phi^H_{21} + \Phi^H_{22}Y_0)(\Phi^H_{11} + \Phi^H_{12}Y_0 )^{-1} - \nonumber
        \\
        \label{eq:Sigma_T_Y0_hybrid}
        &\; (\Phi^H_{11} + \Phi^H_{12}\hat{Y}_0')^{-1} (\Phi^H_{21} + \Phi^H_{22}\hat{Y}_0').
    \end{align}
\end{subequations}
Multiply $\Phi^H_{11} + \Phi^H_{12}\hat{Y}_0'$ from the left and $\Phi^H_{11} + \Phi^H_{12}Y_0$ from the right, we arrive at
\label{eq:Sigma_T_Y0}
\begin{align}
\begin{split}
    & (\Phi^H_{11} + \Phi^H_{12}\hat{Y}_0)'\Sigma_T^{-1} (\Phi^H_{11} + \Phi^H_{12}Y_0) 
    \\
    = & \; (\Phi^H_{11} + \Phi^H_{12}\hat{Y}_0)' (\Phi^H_{21} + \Phi^H_{22}Y_0) 
    \\
    & - (\Phi^H_{21} + \Phi^H_{22}\hat{Y}_0)'(\Phi^H_{11} + \Phi^H_{12}Y_0 ) = Y_0 - \hat{Y}_0, 
\end{split}
\end{align}
where we leveraged the properties of $\Phi^H$ in Lemma \ref{lem:statetransition_hybrid_invertible_E}. From \eqref{eq:boundary_Y0_hatY0}, $Y_0$ and $\hat{Y}_0$ are both symmetric. We can thus let 
\begin{equation}
\label{eq:relation_Y0_Z}
    Y(0) = Z + \frac{1}{2} \Sigma_0^{-1}, \; \hat{Y}(0) = Z - \frac{1}{2} \Sigma_0^{-1},
\end{equation}
for a symmetric $Z$, replace them in \eqref{eq:Sigma_T_Y0_hybrid}, and get
\begin{equation}
\label{eq:Phi_Z_equation}
    (\Phi^H_{11} + \Phi^H_{12} (Z - \frac{\Sigma_0^{-1}}{2} ) )'\Sigma_T^{-1} (\Phi^H_{11} + \Phi^H_{12}(Z + \frac{\Sigma_0^{-1}}{2} )) = \Sigma_0^{-1}.
\end{equation}
Simplifying the above equation, we arrive at 
\begin{equation*}
    \begin{split}
        \!\! Z_{\pm} =& \pm \Sigma_0^{-\frac{1}{2}} (\frac{I}{4} + \Sigma_0^{\frac{1}{2}} (\Phi^H_{12})^{-1} \Sigma_T ((\Phi^H_{12})')^{-1} \Sigma_0^{\frac{1}{2}})^{\frac{1}{2}} \Sigma_0^{-\frac{1}{2}}
        \\
        &-(\Phi^H_{12})^{-1}\Phi^H_{11}. 
    \end{split}
\end{equation*}
By the same arguments as in (\cite{CheGeoPav17a}, Theorem 1) using the monotonicity of $\Phi^H_{12}(t, 0)^{-1}\Phi^H_{11}(t,0)$, $Z_{+}$ has a finite escape time and only $Z_-$ is the meaningful solution to \eqref{eq:Phi_Z_equation}. Plug $Z_{-}$ into \eqref{eq:relation_Y0_Z}, we obtain the initial conditions \eqref{eq:main_result_Pi0_H0}.
\end{proof}
\end{theorem}

\section{Results for Rectangular Jump Dynamics}
\label{sec:E_rectangular}
This section shows the results for the general jump dynamics $\Xi$. In this case, the definition of $\Phi^\Xi$, $H_j(t^+)$ in \eqref{eq:defn_Phi_minus_invertible} and \eqref{eq:hybrid_transition_H} are invalid, and the approach to solve Problem \ref{prob:overall_problem} is thus different from solving Problem \ref{prob:prob_invertible_E}. 

\subsection{The Optimal Distribution of Schr\"odinger Bridge}
Consider first the problem of minimizing \eqref{eq:obj_measure_1} over the path measure $\hat\mP$ induced by the linear process \eqref{eq:dyn_smooth_j} \textit{without} marginal covariance constraints, the unconstrained problem is thus a convex optimization. In mode $I_j$, it is
\begin{equation}
\label{eq:KL_minimization_unconstrained}
     \min_{\hat\mP_j} \langle \epsilon \log d\hat\mP_j - \epsilon d\mQ_j + \mathcal{L}_j , d\hat\mP_j\rangle; \; {\rm s.t. } \langle 1, d\hat\mP_j\rangle = 1.
\end{equation}
Introducing a multiplier $\xi_j$ gives the Lagrangian 
\[
L_{\mathcal{J}_j} = \langle \mathcal{L}_j + \epsilon \log d \hat\mP_j - \epsilon \log d\mQ_j + \xi_j , d \hat\mP_j \rangle - \xi_j.
\]
The necessary condition for the optimal $d\hat\mP^*$ is 
\[
0 = \delta L_{\mathcal{J}_j} |_{d \hat\mP_j^*} = \mathcal{L}_j + \epsilon \left(\log d \hat\mP_j^* - \log d \mQ_j\right) + \xi_j + 1,
\]
leading to $d\hat\mP_j^* \propto d\mQ_j \exp \left( - \mathcal{L}_j / \epsilon \right)$. After normalization, we see that \[
\hat\mP_j^* = \mP_j^*.
\]
We now construct $\mP^*$ explicitly from $\hat\mP^*$ by solving an unconstrained optimal control problem in the following Lemma \ref{lem:closed_form_hatP_star}.
\begin{lemma}
\label{lem:closed_form_hatP_star}
     The objective function of covariance steering problem in mode $I_j$ is equivalent to 
     \[
     \mathcal{J}_j = {\rm KL} ( \mP_j \parallel \mP_j^*),
     \]
     where $\mP_j^*$ is induced by the solution to the LQG problem
    \begin{equation}
    \label{eq:unconstrained_LQG}
        \begin{split}
            &\min_{\hat{u}_j(t)} \;\; \mE \left\{ \int_{t_0^j}^{t_f^j} \lVert \hat{u}_j(t) \rVert^2 + X_j(t)'Q_j(t)X_j(t) dt \right\}
        \end{split}
    \end{equation}
    for systems \eqref{eq:dyn_smooth_j} without marginal constraints.
    
    \begin{proof}
        By Girsanov's theorem \eqref{eq:Girsanov} and the definition of $\hat\mP$, the unconstrained KL-minimization \eqref{eq:KL_minimization_unconstrained} is equivalent to the standard LQG problem \eqref{eq:unconstrained_LQG}. The measure $\hat\mP^*$ is the one induced by the solution to \eqref{eq:unconstrained_LQG}. On the other hand, the covariance steering with quadratic state cost is equivalent to minimizing ${\rm KL}(\mP_j \parallel \mP_j^*)$, and we have shown that $\hat\mP_j^* = \mP_j^*$. 
    \end{proof}
\end{lemma}

The optimal control to the standard LQG problem \eqref{eq:unconstrained_LQG} is 
\[
\hat{u}_j^*(t) = -B_j(t)'\hat\Pi_j(t)X_j(t),
\]
where $\hat\Pi_j(t)$ is the solution to the Riccati equation \eqref{eq:riccati_smooth} with $\hat\Pi_j(t_f^j) = 0$. Denote the resulting closed-loop system as 
\[
\hat{A}_j(t) \triangleq A_j(t)-B_j(t)B_j(t)'\hat\Pi_j(t),
\]
the measure $\mP_j^* (\hat\mP_j^*)$ is thus induced by the process 
\begin{equation}
\label{eq:smooth_dyn_closeloop_unconstrained}
    dX_j(t) = \hat{A}_j(t) X_j(t) d t + \sqrt{\epsilon} B_j(t) dW_j(t).
\end{equation}
We denote the state transition matrix associated with the system $\hat{A}_j(t)$ as $\Phi^{\hat{A}_j}(t,s)$, satisfying $\partial \Phi^{\hat{A}_j}(t,s)/\partial t = A_j(t) \Phi^{\hat{A}_j}(t,s)$ with $\Phi^{\hat{A}_j}(s,s)=I$. Let $\Phi_{\hat{A}_j} \triangleq \Phi^{\hat{A}_j}(t_f^j,t^j_0)$. 

We are ready to reformulate the overall H-CS problem into a convex optimization.

\subsection{A Finite Dimensional Convex Program for H-CS}

In mode $I_j$, denote the joint marginal distribution between the initial and terminal time state under path measure $\mP_j$ by $\mP_{(t_0^j,t_f^j)}$, and denote the same joint marginals for measure $\mP_j^*$ as $\mP_{(t_0^j,t_f^j)}^*$. The objective \eqref{eq:obj_measure_2} can be decomposed into the sum of the relative entropy between the joint initial-terminal time marginals and the relative entropy between the path measures conditioned on these two joint marginals as
\begin{equation}
\label{eq:relative_entropy_marginal}
\begin{split}
    \mathcal{J}_j 
= &\;  {\rm KL}\left( \mP_{(t_0^j,t_f^j)} \parallel \mP_{(t_0^j,t_f^j)}^* \right) 
\\
& + {\rm KL}\left( \mP(\cdot | (t_0^j,t_f^j)) \parallel \mP^*(\cdot | (t_0^j,t_f^j)) \right).
\end{split}
\end{equation}
The relative entropy between the conditional path measures can be chosen to be zero \cite{CheGeoPav15a}, reducing \eqref{eq:Girsanov} into the finite-dimensional problem \eqref{eq:relative_entropy_marginal}. i.e.,
\begin{equation}
\label{eq:obj_finite_dim}
    \mathcal{J}_j = {\rm KL}\left( \mP_{(t_0^j,t_f^j)} \parallel \mP_{(t_0^j,t_f^j)}^* \right).
\end{equation}
For linear SDEs, both $\mP_{(t_0^j,t_f^j)}$ and $\mP_{(t_0^j,t_f^j)}^*$ are Gaussian distributions. Denote the covariance of $\mP_{(t^j_0,t_f^j)}$ as $\Sigma_{t^j_0,t_f^j}$, and denote the covariance of $\mP_{(t_0^j,t_f^j)}^*$ as $\Sigma^{*}_{t^j_0, t_f^j}$. From the solution to linear SDE \eqref{eq:smooth_dyn_closeloop_unconstrained}, the terminal-time state follows
\[
\hat{X}_j(T) = \Phi_{\hat{A}_j}X_0 + \sqrt{\epsilon}\int_{t_0^j}^{t_f^j} \Phi^{\hat{A}_j}(\tau, t^j_0) B_j(\tau) dW_j(\tau),
\]
and the covariance of the initial-terminal joint Gaussian is  
\begin{equation*}
% \label{eq:Sigma_u0_0T}
     \Sigma^{*}_{0,t^-} \!=\! \begin{bmatrix}
         \Sigma_0 & \; \Sigma_0 \Phi_{\hat{A}_1}'
         \\
         \Phi_{\hat{A}_1} \Sigma_0 & \; \Phi_{\hat{A}_1} \Sigma_0 \Phi_{\hat{A}_1}' + \epsilon S_1
     \end{bmatrix}
\end{equation*}
which is a fixed matrix. In $[t^+, T]$, we have 
\begin{equation}
\label{eq:Sigma_u0_0T_A2}
     \Sigma^{*}_{t^+,T} \!=\! \begin{bmatrix}
         \Sigma^+ & \; \Sigma^+ \Phi_{\hat{A}_2}'
         \\
         \Phi_{\hat{A}_2} \Sigma^+ &  \;\Phi_{\hat{A}_2} \Sigma^+ \Phi_{\hat{A}_2}' + \epsilon S_2
     \end{bmatrix}.
\end{equation}
In the above equations, for $j\in\{1,2\}$,
\begin{equation}
\label{eq:Grammian}
\!\!\!\!\!\!S_j \triangleq \int_{t^j_0}^{t_f^j} \Phi^{\hat{A}_j}(\tau,t_0^j) B_j(\tau) B_j'(\tau) (\Phi^{\hat{A}_j}(\tau,t_0^j))' d\tau.
\end{equation}
By the definition of $\mP_{(t^j_0,t_f^j)}$, they have designated initial and terminal covariances. We can thus construct
\begin{equation}
\label{eq:Sigma_u_0T}
\Sigma_{0,t^-} = \begin{bmatrix}
    \Sigma_0 & W_1'
    \\
    W_1 & \Sigma^-
\end{bmatrix}, \; \Sigma_{t^+,T} = \begin{bmatrix}
    \Sigma^+ & W_2'
    \\
    W_2 & \Sigma_T
\end{bmatrix},
\end{equation}
where $\Sigma^-, \Sigma^+, W_1, W_2$ are the optimization variables. The relative entropy $\eqref{eq:obj_finite_dim}$ between two Gaussians is convex and has a closed-form expression in their covariances
\begin{equation}
\label{eq:convex_obj_1}
    \begin{split}
    \mathcal{J}_j = & \; {\rm KL}\left( \mP_{(t_0^j,t_f^j)} \parallel \mP_{(t_0^j,t_f^j)}^* \right) = \!\frac{1}{2} \left(\log\det \Sigma^{*}_{t_0^j,t_f^j} \right.
    \\
    \!\!\!\! & \left. - \log\det \Sigma_{t^j_0,t_f^j} + \tr(\Sigma^{*-1}_{t^j_0,t_f^j}\Sigma_{t^j_0,t_f^j}) \right).
    \end{split}
\end{equation}
We have seen in \eqref{eq:J_decompose} that the objective function for H-CS is the sum of $\mathcal{J}_{j\in\{1,2\}}$ in the above distributional control formulation \eqref{eq:obj_finite_dim}. The additional constraint at the jump event \eqref{eq:cov_constraint_j} for the mean is naturally satisfied by the mean trajectory. In the distribution dual, this constraint becomes the constraint for the covariances \eqref{eq:covariance_evolution_hybrid}. It is exactly 
\begin{equation}
\label{eq:constraint_hat_E}
\Sigma^+ = \Xi \Sigma^- \Xi'.
\end{equation}
The H-CS problem is thus equivalent to minimizing \eqref{eq:convex_obj_1} subject to the constraint \eqref{eq:constraint_hat_E}.

{\em Singularity issues.}
The general saltation $\Xi$ poses issues when covariance at the post-event time instance, $\Sigma^+$, is singular. A notable example is when the system jumps from a lower dimensional state space to a higher dimensional one, i.e., $n_2 > n_1$. This singularity leads to infinite objective functions in \eqref{eq:convex_obj_1}. We show in the sequel that this issue can be avoided in the proposed optimization formulation. 

\subsection{An SDP formulation}
Define a new positive scalar variable $\eta \in \mR$, and a new matrix variable 
\[
\Sigma^+_{\eta} \triangleq \Sigma^+ + \eta I,\; \Sigma^+ = \lim_{\eta \to 0} \Sigma^+_{\eta}.
\] 
In this way, $\Sigma^+_{\eta}$ can be chosen to be non-singular. We replace $\Sigma^+$ with $\Sigma^+_{\eta}$ in all the above expressions \eqref{eq:Sigma_u0_0T_A2}, \eqref{eq:Sigma_u_0T} to continue the following analysis, and let $\eta \to 0$ at the end of the analysis to recover the case for singular $\Sigma^+$. 

Next, we write the objective \eqref{eq:convex_obj_1} as a function explicitly in the variables $W_1, W_2, \Sigma^-, \Sigma^+_\eta$ in the following Lemma \ref{lem:convex_obj}.
\begin{lemma}
\label{lem:convex_obj}
    The objective \eqref{eq:convex_obj_1} is equivalently
    \begin{align}
    \label{eq:objective_convex_2}
    \begin{split}
         & \!\!\!\!\! \mathcal{J}_H \propto \frac{1}{\epsilon}\tr ( S_1^{-1}\Sigma^- ) - \log\det(\Sigma_T - W_2 (\hat\Sigma^{+}_\eta)^{-1} W'_2)
        \\
        & \! - \frac{2}{\epsilon} \tr (\Phi_{\hat{A}_1}' S_1^{-1} W_1 ) + \frac{1}{\epsilon}\tr ( \Phi_{\hat{A}_2}' S_2^{-1} \Phi_{\hat{A}_2} \Sigma^+_\eta )
        \\
        & \! - \log\det ( \Sigma^- - W_1\Sigma_0^{-1}W_1') - \frac{2}{\epsilon} \tr (\Phi_{\hat{A}_2}' S_2^{-1} W_2 )
    \end{split}
    \end{align}
    which is convex in all the variables $W_1, W_2, \Sigma^-, \Sigma^+$. 
    
    \begin{proof}
        See Appendix \ref{sec:proof_lemma_convex_obj}.
    \end{proof}
\end{lemma} 

We notice that the objective \eqref{eq:objective_convex_2} still has a term $(\Sigma_{\eta}^+)^{-1}$. Next, we show in the following Theorem \ref{thm:main_generalE} that this convex problem is a Semi-Definite Programming (SDP), which avoids explicit matrix inversions. 
\begin{theorem}[Main results for general jump dynamics]
    \label{thm:main_generalE}
    The main problem \eqref{eq:formulation_main} with a general hybrid transition $\Xi$ is equivalent to the following Semi-Definite Program
    \begin{subequations}
\label{eq:convex_formulation_sdp}
    \begin{align}
        & \min_{W_1, W_2, \Sigma^-, \Sigma^+, Y_1, Y_2}  \mathcal{J}_{\rm\; SDP}
        \\
         {\rm s.t.} \quad & \Sigma^- \succ 0, \; Y_1 \succ 0, \; \Sigma^+ \succcurlyeq 0, \; Y_2 \succcurlyeq 0,
        \\
        & \Sigma^+ = \Xi \Sigma^- \Xi',
        \\
        & \begin{bmatrix}
            \Sigma^+ & W_2'
            \\
            W_2 & \Sigma_T - Y_2
        \end{bmatrix} \succcurlyeq 0,
        \label{eq:sdp_d}
    \end{align}
\end{subequations} 
where $\mathcal{J}_{\rm\; SDP}$ equals 
\begin{align*}
    \begin{split}
        & \frac{1}{\epsilon}\tr ( S_1^{-1}\Sigma^- ) - \frac{2}{\epsilon} \tr (\Phi_{A_2}' S_2^{-1} W_2 ) - \frac{2}{\epsilon} \tr (\Phi_{A_1}' S_1^{-1} W_1 ) 
        \\
        & \;\; + \frac{1}{\epsilon}\tr ( \Phi_{A_2}' S_2^{-1} \Phi_{A_2} \Sigma^+) - \log\det ( Y_1) - \log\det(Y_2).
    \end{split}
\end{align*}
\begin{proof}
    For the objective function in the expression \eqref{eq:objective_convex_2}, define a matrix slack variable $Y_2\in \mR^{n_2\times n_2}$. By the properties of Schur complement, the following part of the objective  
    \[
    \min_{W_2, \Sigma^+}  - \log\det(\Sigma_T - W_2 (\Sigma^{+}_\eta)^{-1} W'_2);  \;\;  {\rm s.t.} \quad \Sigma^+_\eta \succ 0
    \]
is equivalent to 
\[
\min_{W_2, \Sigma^+_\eta, Y_2}  - \log\det(Y_2);  \quad  {\rm s.t.} \quad \begin{bmatrix}
    \Sigma^+_\eta & W_2'
    \\
    W_2 & \Sigma_T - Y_2
\end{bmatrix} \succcurlyeq 0.
\]
In this way, the term $\Sigma^+_\eta$ appears as a constraint, which can include the singular case. Notice that the above equivalence is true for any $\eta$. Now, to recover the original objective, we let $\eta \to 0$, recovering $\Sigma^+$ in \eqref{eq:sdp_d}.

Similarly, define $Y_1 \triangleq \begin{bmatrix}
    \Sigma_0 & W_1'
    \\
    W_1 &\Sigma^- 
\end{bmatrix} \in \mR^{2 n_1\times 2 n_1}$, the part of the objective 
\begin{equation*}
\label{eq:formulation_convex_5}
        \min_{W_1, \Sigma^-}  - \log\det(\Sigma^- - W_1 \Sigma_0^{-1} W'_1);  \quad  {\rm s.t.} \;\; \Sigma^- \succ 0
\end{equation*}
is then equivalent to 
\[
\min_{\Sigma^-, Y_1} -\log\det(Y_1); \quad {\rm s.t.} \;\; Y_1 \succ 0.
\]
Combining these parts, we arrive at the objective $\mathcal{J}_{SDP}$ in the Theorem, and also the SDP formulation \eqref{eq:convex_formulation_sdp}.
\end{proof}
\end{theorem} 

\begin{remark}
    We note that the above formulation is convex. It solves the general H-CS problem efficiently since the complexity increases only linearly with the number of hybrid transitions instead of the number of time discretizations.
\end{remark}

The solved $\Sigma^-$ and $\Sigma^+$ are optimal for the overall problem with hybrid transitions. After obtaining these covariances, we can recover the optimal controllers by \eqref{eq:initial_a} or \eqref{eq:terminal_a} in the two time windows.

%% file: experiments.tex
\section{Numerical Example}
\label{sec:experiments}
\subsection{Bouncing Ball with Elastic Impact}
We validate the proposed methods with 1D bouncing ball dynamics under elastic impact. The system has $2$ modes, $\mathcal{I} = \{I_1, I_2\}$ where the domains are defined as $D_1 \triangleq \{ (z,\dot z) | \dot z < 0 \}, \; D_2 \triangleq \{ (z,\dot z) | \dot z \geq 0 \}.$ The state in the $2$ modes has the same physical meaning and consists of the vertical position and velocity $X^0=X^1\triangleq[z, \dot z]'$. The control input is a vertical force. 
The system has the same flow 
\[
    \!\!F_j \triangleq dX^j_t = \left(\begin{bmatrix}
        0 & 1
        \\
        0 & 0
    \end{bmatrix}X^j_t + \begin{bmatrix}
        0 \\ \frac{u^j_t-mg}{m}
    \end{bmatrix} \right) dt + \sqrt{\epsilon}\begin{bmatrix}
        0 \\ \frac{1}{m}
    \end{bmatrix} dW^j_t.
\]
The guard functions are defined as $g_{12}(z,\dot z) \triangleq z$, and $g_{21} \triangleq \dot z$. The reset map $R_{21}$ is the identity map, and 
\[
    \begin{bmatrix}
     z(t^+) 
     \\
     \dot z(t^+)
    \end{bmatrix}
    =
    R_{12}(z(t^-), \dot z(t^-)) \coloneqq 
    \begin{bmatrix}
    1 & 0 
    \\
    0 & -e_2
    \end{bmatrix}
    \begin{bmatrix}
     z(t^-) 
     \\
     \dot z(t^-)
    \end{bmatrix}
\]
where $e_2$ represents the elastic loss in position and velocity at the contact. In our case, we choose $e_2 = 0.6$. The bounding ball dynamics corresponds to the invertible reset map case introduced in Section \ref{sec:E_invertible}. The start and goal mean states are $x_{S} = [5, 1.5]'; \; x_{G} = [2.5, 0]'$, both in mode $I_2$. The initial and the target covariances are
\[
    \Sigma_0 = 0.2\mathrm{I}, \; \Sigma_T = 0.05\mathrm{I}.
\]
All experiments use $\Delta t=0.0015, \epsilon=0.5$. We use the H-iLQR controller \cite{kong2021ilqr, kong2023hybrid} to obtain the mean trajectory by setting a quadratic terminal loss $\Psi_T(X_T) \coloneqq \frac{1}{2} \lVert X_T \rVert^2_{Q_T}$ where $Q_T=25 \mathrm{I}$, and then we update the controller using the optimal covariance steering feedback \eqref{eq:optimal_feedback}.

\begin{figure}[th]
    \centering
    \begin{subfigure}{0.65\linewidth}
    \centering
        \includegraphics[width=\linewidth]{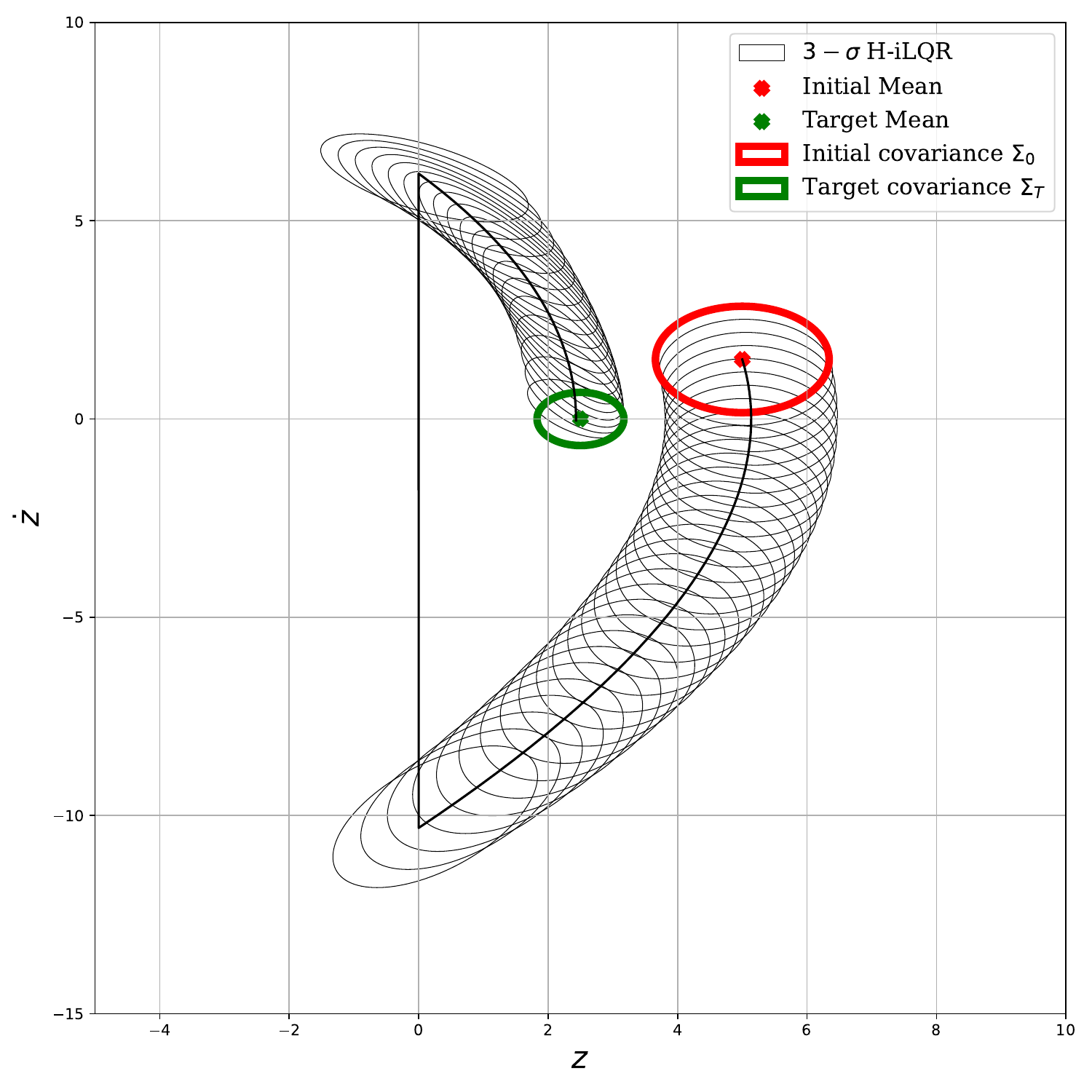}
    \caption{Controlled Covariance under H-iLQR.}
    \end{subfigure}
    \begin{subfigure}{0.65\linewidth}
    \centering
        \includegraphics[width=\linewidth]{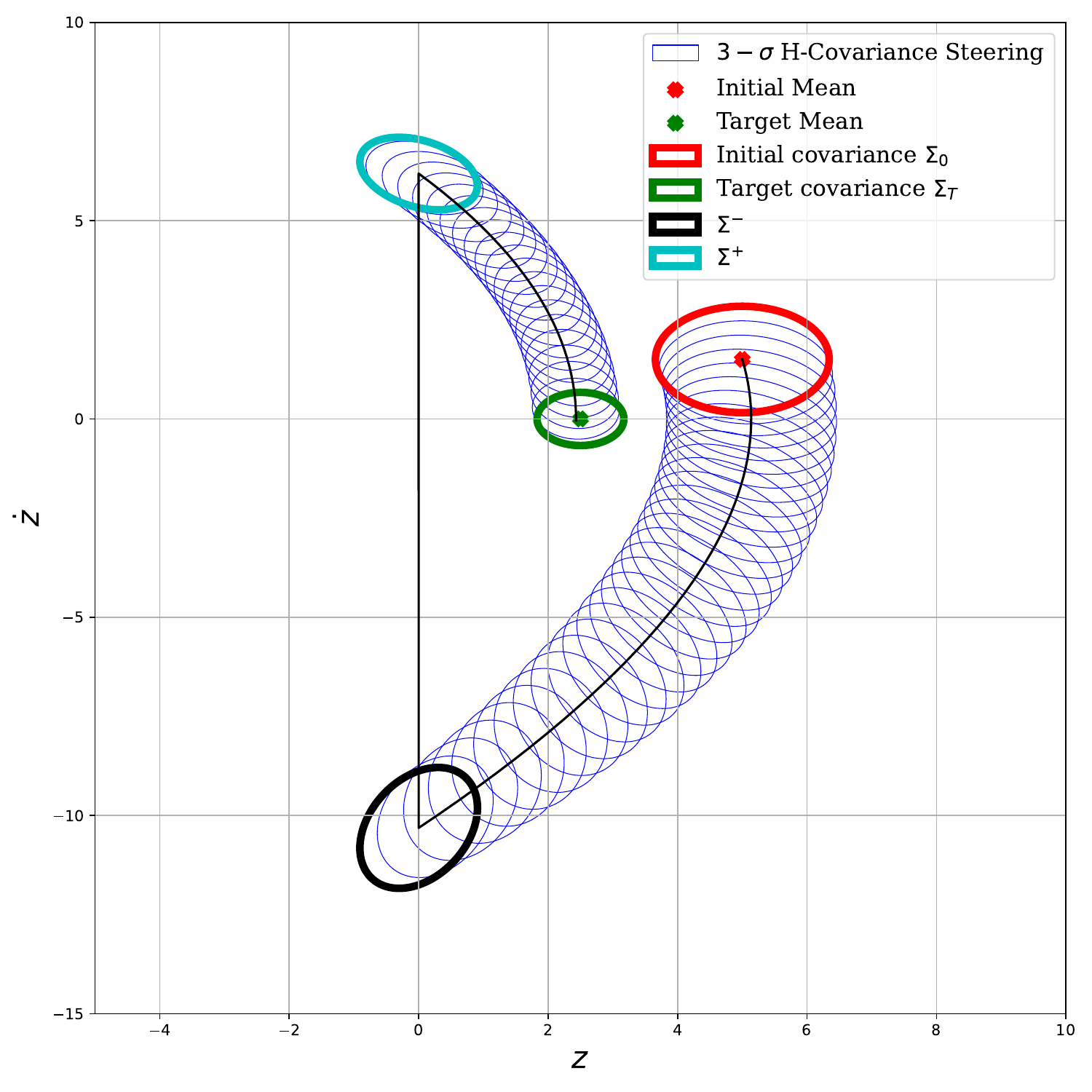}
    \caption{Controlled Covariance under H-CS.}
    \end{subfigure}
    \caption{Covariance steering for a bouncing ball dynamics with elastic impacts. The H-CS controller guarantees the terminal covariance constraint.}
    \label{fig:bouncing}
\end{figure}

\begin{figure}[th]
    \centering
\includegraphics[width=1.0\linewidth]{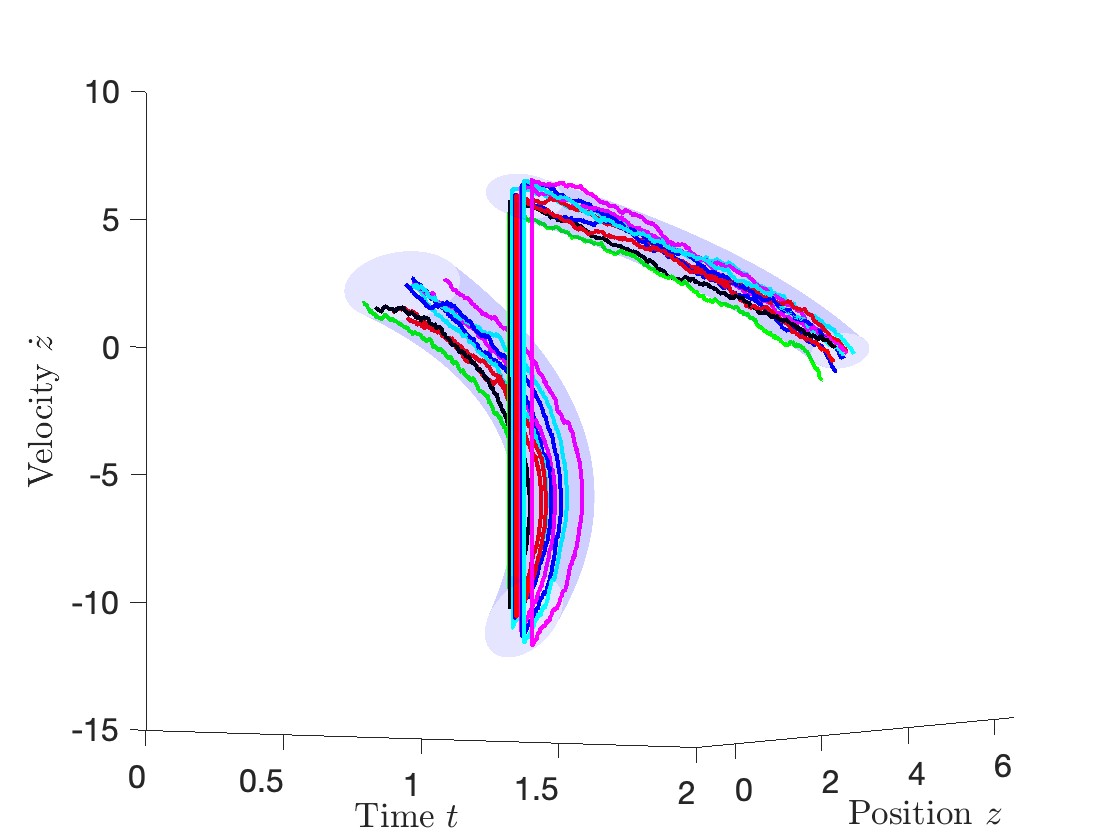}
    \caption{Controlled covariance tube and sampled trajectories for the bouncing ball dynamics.}
    \label{fig:bouncing_samples}
\end{figure}

The results are shown in Fig. \ref{fig:bouncing}. In Fig. \ref{fig:bouncing}, the mean is obtained from the H-iLQR controller, and the covariance control is achieved by integrating the hybrid Riccati equations \eqref{eq:riccati_smooth}, \eqref{eq:riccati_hybrid} starting from the initial conditions \eqref{eq:main_result_Pi0_H0}. We compute the covariance propagation under the H-iLQR controller, plotted in black for comparison. H-CS controller steered the state covariance to the target terminal one. Fig. \ref{fig:bouncing_samples} plots the controlled $3-\sigma$ covariance tube and sampled trajectories under the same controller. We observe that the bouncing timing is affected by the uncertainties. The distribution arrives at the target covariance at the terminal time. We note here that the H-iLQR controller can steer the state covariance to one that smaller than the target one at the terminal time, but it requires additional penalties added to the terminal costs, which corresponds to additional control energy costs. The H-CS gives the exact cost design to steer the state covariance to a designated one, corresponding to the minimum amount of control energy needed to accomplish the uncertainty control task. 

We use the SDP formulation \eqref{eq:convex_formulation_sdp} to solve the same problem, and compare $\Sigma^-, \Sigma^+$ with the analytical solution in the Table \ref{tab:cov_comparison_bouncing}. Results show that the optimization recovers the analytical solutions. Both methods verify the conditions \eqref{eq:covariance_evolution_hybrid}.
\begin{table}[ht]
    \centering
    \caption{The analytical and the numerical solutions for the pre-event and post-event covariances.}
    \label{tab:cov_comparison_bouncing}
    \begin{tabular}{|*3{>{\renewcommand{\arraystretch}{1}}c|}}
    \hline
         &  Analytical (Theorem \ref{thm:main}) & Numerical (Theorem \ref{thm:main_generalE})
         \\
         \hline
       $\Sigma^-$  & 
       $\begin{bmatrix}
       0.1651 & 0.08517
       \\
        0.0852 & 0.2640
       \end{bmatrix}$
       & 
       $\begin{bmatrix}
        0.1651 & 0.0855
        \\
        0.0855 & 0.2640
       \end{bmatrix}$
         \\
         \hline
         $\Sigma^+$  
         & 
         $\begin{bmatrix}
            0.1651 & -0.0513 
            \\
            -0.0513 &  0.0950
         \end{bmatrix}  $
         & 
         $\begin{bmatrix}
              0.1651 & -0.0511 
              \\
            -0.0511 &  0.0950
         \end{bmatrix}$
         \\
         \hline
    \end{tabular}
\end{table}

\subsection{Linearized Spring-Loaded Inverted Pendulum (SLIP): Jumping to a higher dimensional state space.}
\label{sec:experiments_slip}
We validate our algorithm on a linearized classical Spring-Loaded Inverted Pendulum (SLIP) dynamics. We linearize the SLIP model around the mean trajectory.  

SLIP consists of a body modeled by a point mass $m$ and a leg modeled by a massless spring with coefficient $k$ and natural length $r_0$. We use $(p_x, p_z)$ to represent the horizontal and vertical positions of the body, and $(v_x, v_z)$ are their velocities, respectively. $(\theta, \dot \theta)$ are the angle and angular velocity between the leg and the horizontal ground, and $(r, \dot r)$ are the leg length and its changing rate. The variables $[p_x, v_x, p_z, v_z, \theta, \dot \theta, r, \dot r]$ fully describe the system's states.  

The SLIP system has two modes, denoted as $\mathcal{I} = \{I_1, I_2\}$ where the domains are defined as $D_1 \triangleq \{ X_1(t) | r - r_0 < 0  \}; \; D_2 \triangleq \{ X_2(t) | p_z - r_0 \sin{\theta} \geq 0 \}.$ $I_1$ is also known as the \textit{stance} mode, and mode $I_2$ is known as the \textit{flight} mode. In the stance mode, the positions $p_x$ and $p_z$ are constrained by the spring-loaded leg, leading to different state spaces in the two modes. We use a polar coordinate system for the stance mode. 
In the stance mode, the state is represented by $X_1(t) \triangleq [\theta, \dot \theta, r, \dot r]$, and the stance smooth flow is
\begin{align*}
    \begin{split}
dX_1 = &\left(\begin{bmatrix}
        \dot \theta\\
        \frac{-2\dot \theta\times \dot r - g\times \cos\theta}{r}\\
        \dot r\\
        \frac{k(r0-r)}{m} - g\sin\theta + \dot \theta^2 r
    \end{bmatrix}
    +
    \begin{bmatrix}
        0 & 0 \\
        0 & 0 \\
        \frac{m}{r^2} & 0  \\
        0 & \frac{k}{m} 
    \end{bmatrix} u_1 \right) d t
    \\
    & + \sqrt{\epsilon} \begin{bmatrix}
        0 & 0 \\
        0 & 0 \\
        \frac{m}{r^2} & 0  \\
        0 & \frac{k}{m} 
    \end{bmatrix} dW_1(t). 
    \end{split}
\end{align*}
In the flight mode, the state $X_2(t) \triangleq [p_x, v_x, p_z, v_z, \theta]$, with the smooth flow 
\begin{equation*}
    d X_2(t) = 
    \begin{bmatrix}
        \dot v_x\\
         0\\
        \dot v_z\\
         -9.81\\
        0
    \end{bmatrix} d t
    +
    \begin{bmatrix}
        0 & 0 & 0\\
        1 & 0 & 0\\
        0 & 0 & 0\\
        0 & 1 & 0\\
        0 & 0 & 1
    \end{bmatrix}
    \left(
    u_2(t) d t + \sqrt{\epsilon}
    dW_2(t)
    \right)  
\end{equation*}
where $u_2(t) \in \mR^3$ is the control input consisting of two acceleration inputs for the body, and an angular velocity input for the leg. The leg is defined by only the state $\theta$ since we assume it to be massless. We assume a double integrator model in the flight mode. 
The reset maps between the two modes are defined as 
\begin{equation*}
\!\!X_2(t^+) = R_{12}(X_1(t^-)) = 
\begin{bmatrix}
    p_{x,T}^- + r_0\cos\theta^- \\
    \dot{r}^-\cos\theta^- - r^-\dot{\theta}^-\sin\theta^- \\
    r_0\sin\theta^-\\
    r_0\dot\theta\cos\theta^- + \dot{r}^-\sin\theta^-\\
    \theta^-
\end{bmatrix},
\end{equation*}
and 
\begin{equation*}
\!\!X_1(t^+) = R_{21}(X_2(t^-)) = 
\begin{bmatrix}
    \theta^- \\
    (p_x^- v_z^- - p_z^- v_x^-) / r_0^2 \\
    r_0\\
    -v_x^-\cos \theta^- + v_z^-\sin\theta^-
\end{bmatrix}
\end{equation*}
where $p_{x,T}$ is the tole $x$-position and is assumed to be unchanged during the stance mode.
        
In our experiment, we choose the dynamics parameters $r_0=1, m=0.5, k=25.0$, and gravity is $-9.81$. Time window is $[0, 0.5]$ with a discretization $\Delta t = 5\times 10^{-5}$. Noise level $\epsilon = 0.0015$. We use the H-iLQR to obtain the mean trajectory, the linearized system, and the Saltation Matrix from the last iteration's backward pass in H-iLQR. 

We set initial state $m_0 \in I_1$ and the target state $m_T \in I_2$ as
\begin{align*}
m_0 &= [1.745, -4.0, 0.5, 0.0],
\\
m_T &= [1.1, 2.25, 1.4, 0.0, \pi/3].
\end{align*}
The resulted mean trajectory is obtained by setting a terminal quadratic cost $Q_T = 2.0 \times \mathrm{I}_5.$ For the covariance control, the initial and terminal covariances are set to be 
\[
\Sigma_0 = 0.002 \mathrm{I}_4, \; \Sigma_T = 0.0003 \mathrm{I}_5.
\]
At the hybrid event time, the system jumps from a lower dimensional state space in mode $I_1$ to a higher dimensional state space in mode $I_2$, leading to a singular post-event covariance matrix. The Saltation matrix we obtain is
\[
\Xi = \Xi_{12} = 
\begin{bmatrix}
    -0.9371 & 0.0   &   0.3492 & 0.0    
    \\
    -2.8513 & -0.9371 & 1.0971 & 0.3492
    \\
    0.3492 & 0.0 &  0.9371 & 0.0  
    \\
    2.3119 & 0.3492 & -0.4088 & 0.9371
    \\
   1.  &    0.0  &    0.3365 & 0.0       
\end{bmatrix}.
\]
The optimized $\Sigma^+$ is 
\[
\Sigma^+ = 
\begin{bmatrix}
  0.0004 & -0.0002 & 0.0   &   0.0   &  -0.0003
  \\
 -0.0002 & 0.0007 & -0.0001 & 0.0001 & 0.0001
 \\
 0.0  &   -0.0001 & 0.0003 & -0.0001 & 0.0002
 \\
 0.0   &   0.0001 & -0.0001 & 0.0007 & -0.0001
 \\
 -0.0003 & 0.0001 & 0.0002 & -0.0001 & 0.0004
\end{bmatrix},
\]
which is indeed singular. The computed $\Sigma^-$ and $\Sigma^+$ satisfy the condition \eqref{eq:covariance_evolution_hybrid}. Using these two covariances in the corresponding smooth covariance steering problems, we obtain the $\Pi(t^-)$ and $\Pi(t^+)$ pair as follows. For the singular covariance $\Sigma^+$, we use the terminal conditions for $\Pi(T)$ in the Remark \ref{thm:remark_terminal_Pi} and integrate the Riccati equation backward to get $\Pi(t^+)$. 

\begin{figure}
    \centering
    \includegraphics[width=\linewidth]{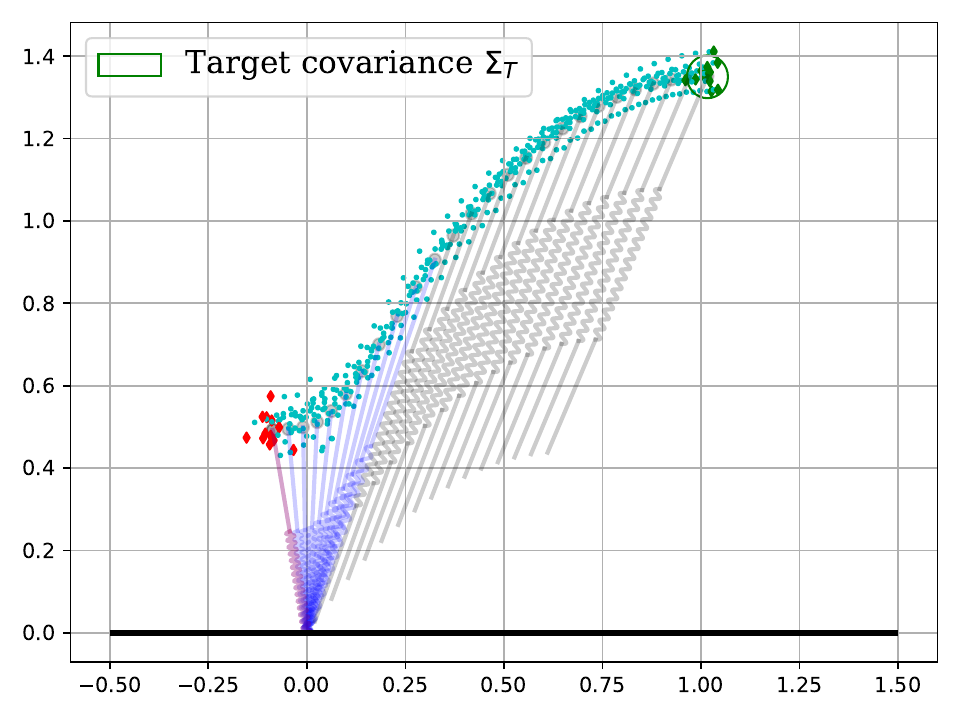}
    \caption{Deterministic nominal trajectory under H-iLQR controller and stochastic trajectories under the H-CS controller for the SLIP model. The transparent legs with springs are the nominal trajectory, and the dots in cyan color are the SLIP body's position trajectories of the samples starting from the initial distribution, marked by red diamond shapes. We draw the $3-\sigma$ boundary of the terminal-time target covariance.}
    \label{fig:samples_SLIP}
\end{figure}

We sample the initial states from the initial Gaussian, and apply the H-CS feedback controller to obtain the samples' state trajectories under uncertainties. Fig. \ref{fig:samples_SLIP} shows the nominal and the controlled stochastic body position trajectory samples. The terminal covariance under the H-iLQR controller, denoted as $\hat\Sigma_T$, is
\[
\hat\Sigma_T = \begin{bmatrix}
    0.0006 & 0.0001 & 0.0001 & 0.0  &  -0.0002
    \\
  0.0001 & 0.0007 & 0.0   &   0.0 &   0.0 
  \\
  0.0001 & 0.0   &   0.0006 & 0.0001  & 0.0002
  \\
 0.0  &  0.0  &  0.0001 & 0.0008 & -0.0001
 \\
 -0.0002 & 0.0  & 0.0002 & -0.0001 & 0.0005
\end{bmatrix},
\]
which does not satisfy the terminal-time covariance constraints.

%% file: conclusion.tex
\section{Conclusion}\label{sec:conclusion}
This work considers the problem of controlling uncertainty for hybrid dynamical systems. We addressed this problem by optimally steering the state covariance of a linear stochastic system subject to hybrid transitions, around a nominal mean trajectory. The nonlinear and instant hybrid transitions can be well-approximated linearly using the Saltation Matrix. For non-singular hybrid transitions, we show that this problem enjoys an analytical solution in closed-form. For the general hybrid transitions, the problem can be re-formulated into a convex optimization with a hybrid equality constraint on the covariance propagation during the hybrid transition event. The optimization is linear in the number of hybrid events and is therefore efficient to solve. Combining the two formulations, we thus showed the complete solution to the covariance steering problem subject to hybrid transitions.  Numerical experiments validated the proposed method. Our method can serve as a feedback component for uncertainty control on top of any given mean trajectory controllers for hybrid dynamics. The future works include iterative covariance control for nonlinear dynamical systems subject to hybrid transitions.

%% file: Appendix.tex
\appendix

\section{Proof of Lemma \ref{lem:lemma_hybrid_Phi_invertible_E}}
\label{sec:appendix_proof_Lemma_hybrid_Phi}
\begin{proof}
For $I_j, j\in\{1,2\}$, a direct time derivatives of \eqref{eq:sol_riccati_XY} for the respective smooth linear systems $(A_j, B_j)$ gives
\begin{equation*}
- \frac{d}{dt} (Y_jX_j^{-1}) = - \dot{Y}X_j^{-1} + Y_jX_j^{-1}\dot{X_j}X_j^{-1}
\end{equation*}
which verifies the Riccati equations \eqref{eq:riccati_smooth} in the respective modes. At the event time $t=t^{-}$, by the definition of $\Phi^{R}(t^+, t^-)$, 
\begin{equation*}
    \begin{bmatrix}
        X(t^+) \\ Y(t^+)
    \end{bmatrix}
    = \begin{bmatrix}
         \Xi & 0
         \\
         0 & (\Xi')^{-1}
        \end{bmatrix}
    \begin{bmatrix}
        X(t^-) \\ Y(t^-)
    \end{bmatrix},
\end{equation*}
and the value function
\begin{equation*}
\begin{split}
    \Pi_1(t^-) &= Y_1(t^{-})X_1^{-1}(t^{-}) = \Xi' Y_2(t^{+})X_2^{-1}(t^{+}) \Xi
    \\
    &= \Xi' \Pi_2(t^+) \Xi
\end{split}
\end{equation*}
verifies \eqref{eq:riccati_hybrid}. Starting from $t^+$ to $t$, the evolution of $\Pi_2(t)$ of the smooth flow governed by $\Phi_2(t)$ has the same time derivation as \eqref{eq:riccati_smooth}. The proof for the evolution of $H(t)$ function follows similar arguments and is neglected.
\end{proof}

\section{Proof of Lemma \ref{lem:state_transition_product}}
\label{sec:appendix_lemma_product}
\begin{proof}
We show the result for equation \eqref{eq:lemma_statetransition_1}. Denote 
    \[
        \bar{\Phi}(t,s) = \left[
    \begin{matrix}
    \bar{\Phi}_{11}(t,s) & \bar{\Phi}_{12}(t,s)\\
    \bar{\Phi}_{21}(t,s) & \bar{\Phi}_{22}(t,s)
    \end{matrix}\right],
    \]
    and 
    \[
    \Tilde{\Phi}(t,s) = \left[
    \begin{matrix}
    \Tilde{\Phi}_{11}(t,s) & \Tilde{\Phi}_{12}(t,s)\\
    \Tilde{\Phi}_{21}(t,s) & \Tilde{\Phi}_{22}(t,s)
    \end{matrix}\right].
    \]
    We directly calculate  
    \begin{align*}
    &  \Phi_{11}(t,s)^{'}\Phi_{22}(t,s) - \Phi_{21}(t,s)^{'}\Phi_{12}(t,s)
    \\
    = &  \;\left(\Tilde{\Phi}_{11}^{'}\bar{\Phi}_{11}^{'} + \Tilde{\Phi}_{21}^{'}\bar{\Phi}_{12}^{'} \right) \left( \bar{\Phi}_{21} \Tilde{\Phi}_{12} +  \bar{\Phi}_{22} \Tilde{\Phi}_{22} \right)
    \\
    - & \;\left(\Tilde{\Phi}_{11}^{'} \bar{\Phi}_{21}^{'} + \Tilde{\Phi}_{21}^{'}\bar{\Phi}_{22}^{'} \right) \left( \bar{\Phi}_{11} \Tilde{\Phi}_{12} + \bar{\Phi}_{12} \Tilde{\Phi}_{22} \right)
    \\
    = & \; \Tilde{\Phi}_{11}^{'} \left( \bar{\Phi}_{11}^{'} \bar{\Phi}_{21} - \bar{\Phi}_{21}^{'} \bar{\Phi}_{11} \right) \Tilde{\Phi}_{12} 
    \\
    +  & \;\Tilde{\Phi}_{11}^{'} \left( \bar{\Phi}_{11}^{'} \bar{\Phi}_{22} - \! \bar{\Phi}_{21}^{'} \bar{\Phi}_{12} \right) \Tilde{\Phi}_{22}
    \\
    + & \;\Tilde{\Phi}_{21}^{'} \left( \bar{\Phi}_{12}^{'} \bar{\Phi}_{21} -  \bar{\Phi}_{22}^{'} \bar{\Phi}_{11} \right) \Tilde{\Phi}_{12} 
    \\
    + & \;\Tilde{\Phi}_{21}^{'} \left( \bar{\Phi}_{12}^{'} \bar{\Phi}_{22} -  \bar{\Phi}_{22}^{'} \bar{\Phi}_{12} \right) \Tilde{\Phi}_{22}
    \\
    = & \; \Tilde{\Phi}_{11}^{'}(t,s) \Tilde{\Phi}_{22}(t,s) - \Tilde{\Phi}_{21}^{'}(t,s) \Tilde{\Phi}_{12}(t,s) = I,
    \end{align*}
    where we used the assumption that both $\Tilde{\Phi}$ and $\bar{\Phi}$ satisfy \eqref{eq:lemma_statetransition}. We thus prove the first equation in \eqref{eq:lemma_statetransition_1}. The other equations in \eqref{eq:lemma_statetransition} can be derived similarly and is omitted.
\end{proof}

\section{Proof of Lemma \ref{lem:statetransition_hybrid_invertible_E}.}
\label{sec:proof_statetransition_hybrid_invertible_E}
\begin{proof}
    The smooth transition kernels $\Phi_1(t^-, 0)$ and $\Phi_2(t, t^+)$ satisfies \eqref{eq:lemma_statetransition} by Lemma \ref{lem:state_transition}. For $t>t^-$, $\Phi^H(t,0)$ equals the product of $\Phi_2(t, t^+)$, $\Phi^{\Xi}(t^+,t^-)$, and $\Phi_1(t^-, 0)$. Given Lemma \ref{lem:state_transition_product}, we only need to prove that $\Phi^{\Xi}$ satisfies \eqref{eq:lemma_statetransition}. A direct computation gives
    \begin{subequations}
        \begin{align*}            (\Phi^{\Xi}_{11})'\Phi^{\Xi}_{22} - (\Phi^{\Xi}_{21})'\Phi^{\Xi}_{12} &= \Xi' (\Xi')^{-1} = I,
            \\
        (\Phi^{\Xi}_{12})'\Phi^{\Xi}_{22} - (\Phi^{\Xi}_{22})'\Phi^{\Xi}_{12} &=(\Phi^{\Xi}_{21})'\Phi^{\Xi}_{11} - (\Phi^{\Xi}_{11})'\Phi^{\Xi}_{21} = 0,
            \\
            \Phi^{\Xi}_{11}(\Phi^{\Xi}_{22})' - \Phi^{\Xi}_{12}(\Phi^{\Xi}_{21})' &= \Xi\Xi^{-1} = I,
            \\
            \Phi^{\Xi}_{12}(\Phi^{\Xi}_{11})' - \Phi^{\Xi}_{11}(\Phi^{\Xi}_{12})' &= \Phi^{\Xi}_{21}(\Phi^{\Xi}_{22})' - \Phi^{\Xi}_{22}(\Phi^{\Xi}_{21})' = 0.
        \end{align*}
\label{eq:lemma_statetransition_invertible_E}
    \end{subequations}
We just showed that $\Phi^{\Xi}(t^+,t^-)$ preserves \eqref{eq:lemma_statetransition_1}. The other equations \eqref{eq:lemma_statetransition_2} - \eqref{eq:lemma_statetransition_6} are straightforward to see by similar calculations. 

In both modes $I_j, j=1, 2$, the entries $\Phi^j_{12}(t^-, s)$ and $\Phi^j_{11}(t^-, s)$ are invertible, and the terms $\Phi^j_{11}(t,s)^{-1} \Phi^j_{12}(t,s)$ are monotonically decreasing functions of $t$, and when $t \searrow 0$, $\Phi^H$ reduces to $\Phi^1$ and has a left limit $0$. At time $t=t^+$, the term
\begin{align*}
    \Phi^H_{12}(t^+, s) &= \Phi^{\Xi}_{11} \Phi^1_{12}(t^-, s) + \Phi^{\Xi}_{12} \Phi^1_{22}(t^-, s) 
    \\
    &= \Xi \Phi^1_{12}(t^-, s)
\end{align*}
is invertible because both $\Xi$ and $\Phi^1_{12}(t^-, s)$ are invertible. Similarly, the term 
\begin{align*}
    \Phi^H_{11}(t^+, s) &= \Phi^{\Xi}_{11} \Phi^1_{11}(t^-, s) + \Phi^{\Xi}_{12} \Phi^1_{12}(t^-, s) 
    \\
    &= \Xi \Phi^1_{11}(t^-, s)
\end{align*}
is invertible. $\Phi^H_{12}(t, s)$ and $\Phi^H_{11}(t, s)$ are thus invertible for all $t$.
Moreover, the term 
\begin{align*}
    \Phi^H_{11}(t^+,0)^{-1} \Phi^H_{12}(t^+,0) &= \Phi^1_{11}(t^-, s)^{-1} \Xi^{-1} \Xi \Phi^1_{12}(t^-, s)
    \\
    &= \Phi^1_{11}(t^-, s)^{-1} \Phi^1_{12}(t^-, s)
\end{align*}
is monotonically decreasing from $t^-$ to $t^+$. We conclude that $\Phi^H_{11}(t,s)^{-1} \Phi^H_{12}(t,s)$ is monotonically decreasing in $[0,T]$.
\end{proof}

\section{Proof of Lemma \ref{lem:convex_obj}}
\label{sec:proof_lemma_convex_obj}
\begin{proof}
    We first calculate the terms $\sum_j (\log\det \Sigma^{*}_{t_0^j,t_f^j} - \log\det \Sigma_{t^j_0,t_f^j} )$. For the fixed process $\mP^*$, in $[0, t^-]$, the term $\Sigma^{*}_{0,t^-}$ is not a function of the variables. In $[t^+, T]$, we have
\[
\log\det \Sigma^{*}_{t^+, T} = \log\det \hat\Sigma^+_{\eta} + \log\det \epsilon S_2 \propto \log\det \hat\Sigma^+_{\eta},
\]
since $\log\det S_2$ is not a variable. For the controlled process $\mP$, in $[0, t^-]$ we have 
\begin{equation}
\label{eq:logdet_Sig_p1}
\!\!\!\! \log\det \Sigma_{0, t^-} = \log\det \Sigma_0 - \log\det\left( \Sigma^{-} - W_1\Sigma_0^{-1}W_1' \right) 
\end{equation}
which is convex in $\Sigma^-$. In time $[t^+, T]$, 
\begin{align*}
& \log\det \Sigma_{t^+, T}  
\\
= & \log\det \Sigma^+_{\eta} + \log\det\left( \Sigma_T - W_2(\Sigma^+_{\eta})^{-1}W_2' \right).
\end{align*}
The terms $\log\det \Sigma^{*}_{t^+,T}$ and $\log\det \Sigma_{t^{+},T}$ both contain $\Sigma^+_{\eta}$, which becomes infinity for singular $\Sigma^+$, and when $\eta \to 0$. However, the difference between them is
\begin{equation}
\label{eq:diff_logdet}
\begin{split}
    &\log\det \Sigma^{*}_{t^+,T} - \log\det \Sigma_{t^{+},T}  
    \\
  \propto &- \log\det \left( \Sigma_T - W_2 (\hat\Sigma^{+}_{\eta})^{-1} W'_2 \right),
\end{split}
\end{equation}
where the term $\log\det \Sigma^+_{\eta}$ is canceled out, and the result is finite and convex in $\Sigma^+_{\eta}$, regardless of the value of $\eta$. 

We next calculate the matrix trace terms in \eqref{eq:convex_obj_1}. In mode $I_2$, applying the block matrix inverse on $\Sigma^{*}_{t^+,T}$ gives
\[
\Sigma^{*-1}_{t^+,T} = \begin{bmatrix}
    (\Sigma^+_{\eta})^{-1} + \frac{1}{\epsilon}\Phi_{\hat{A}_2}' S_2^{-1} \Phi_{\hat{A}_2} & -\frac{1}{\epsilon}\Phi_{\hat{A}_2}'S_2^{-1}
    \\
    -\frac{1}{\epsilon}S_2^{-1}\Phi_{\hat{A}_2} & \frac{1}{\epsilon}S_2^{-1}
\end{bmatrix},
\]
leading the term $\tr(\Sigma^{*-1}_{t^+,T} \Sigma_{t^+,T})$ in \eqref{eq:convex_obj_1} to be 
\[
\begin{split}
    \tr(\Sigma^{*-1}_{t^+,T} \Sigma_{t^+,T}) &= \tr \left(
    ((\Sigma^+_{\eta})^{-1} + \frac{1}{\epsilon}\Phi_{\hat{A}_2}'S_2^{-1}\Phi_{\hat{A}_2} )\Sigma^+_{\eta}  \right.
    \\
    & \!\!\!\!\!\!\!\!\!\!\!\!\! \left. -\frac{1}{\epsilon}\Phi_{\hat{A}_2}'S_2^{-1}W_2 - \frac{1}{\epsilon}S_2^{-1}\Phi_{\hat{A}_2}W_2' + \frac{1}{\epsilon}S_2^{-1}\Sigma_T
    \right)
    \\
    & \!\!\!\!\!\!\!\!\!\!\!\!\!\!\!\!\!\!\!\!\!\!\!\! \propto \tr \left(\frac{1}{\epsilon}\Phi_{\hat{A}_2}'S_2^{-1}\Phi_{\hat{A}_2} \Sigma^+_{\eta}  - \frac{2}{\epsilon}\Phi_{\hat{A}_2}'S_2^{-1}W_2\right).
\end{split}
\]
The term $(\Sigma^+_{\eta})^{-1}$ is canceled out in this part of the objective function. After a similar calculation for $\tr(\Sigma^{*-1}_{0,t^-} \Sigma_{0,t^-})$, and combining \eqref{eq:logdet_Sig_p1} and \eqref{eq:diff_logdet}, we obtain the result.
\end{proof}